\documentclass[12pt,oneside]{amsart}

\pdfoutput=1

\usepackage[section]{placeins}
\usepackage[margin=1in]{geometry}
\usepackage{mathabx}
\usepackage{xcolor}
\usepackage{float}
\usepackage{graphicx}
\usepackage{amssymb}
\usepackage{tikz}
\usepackage{placeins}
\usepackage{amsmath,amssymb}
\usetikzlibrary{arrows.meta}
\usepackage{latexsym}
\usepackage{epsfig}
\usepackage{enumitem}
\usepackage{multicol}
\usepackage{wasysym}
\usepackage{amsmath}
\usepackage{amsthm}
\usepackage{amscd}
\usepackage{ulem}
\usepackage{soul}
\usepackage{multirow}
\usepackage{rotating}
\usepackage{pdflscape}
\usepackage{float}
\usepackage{dsfont}
\usepackage{comment}

\usepackage[colorlinks,citecolor=red,pagebackref,hypertexnames=false]{hyperref}
\hypersetup{
  colorlinks = true,
  citecolor=red,
  linkcolor  = red
}

\setlist{itemsep=.06125in}

\restylefloat{table}
\numberwithin{equation}{section}

\theoremstyle{plain}
\newtheorem{theorem}{Theorem}[section]
\newtheorem{lemma}[theorem]{Lemma}
\newtheorem{corollary}[theorem]{Corollary}
\newtheorem{proposition}[theorem]{Proposition}

\theoremstyle{definition}
\newtheorem{definition}[theorem]{Definition}

\newtheorem{example}[theorem]{Example}

\theoremstyle{remark}
\newtheorem{remark}[theorem]{Remark}

\def\R{\mathbb R}

\def\S{\mathbb{S}}

\date{\today}

\author{Shantanu Deodhar}
\address{Department of Mathematics, University of Rochester, Rochester, NY, USA}
\email{sdeodhar@ur.rochester.edu}
\subjclass[2020]{Primary 42B10; Secondary 42B20, 28A75}
\keywords{$L^{p}$ integrability, Spectral Synthesis, Fourier Restriction, Quadratic Manifolds, Fourier Ratio}

\begin{document}
\title{$L^{p}$ integrability and spectral complexity of  manifolds with higher co-dimension}
\begin{abstract}
We investigate the relationship between the $L^p$ integrability
of functions and the curvature geometry of their Fourier support.
We prove that if $M\subset\R^n$ is a $C^n$ manifold of dimension
$d$ and $n>2d$, then every $u\in L^p(\R^n)$ with
$\operatorname{supp}\widehat u\subseteq M$ vanishes whenever
\[
2\le p\le\frac{6n-4}{3d-1}.
\]
This strictly improves the classical Agranovsky-Narayanan threshold $2n/d$ throughout
the range $n>2d$. More generally, we obtain a family of
dimension-dependent improvements in higher codimension.
As a consequence, we obtain a strict quantitative improvement
in the necessary conditions for Fourier restriction when
$n>2d$, resolving a conjecture of Mockenhaupt
\cite{Mockenhaupt96}. This further establishes the nonexistence
of quadratic manifolds satisfying the Christ--Mockenhaupt
condition when $n>2d$, answering the question of \cite{CLMP}.

For quadratic manifolds, we obtain spectral synthesis bounds
from joint sublevel estimates for the eigenvalues of the normal
Hessian pencil. In codimension two, we determine the exact
critical exponent from the ranks and eigenvalue vanishing orders,
including nonexistence at the critical endpoint. We also obtain
integrability bounds for manifolds with constant relative nullity,
extending Guo's hypersurface results \cite{Guo93} to higher
codimension. Finally, we use the Fourier ratio as a measure
of the geometric complexity of individual strata of the
characteristic set to establish $L^p$ uniqueness criteria
for homogeneous constant-coefficient PDEs.
\end{abstract}
\maketitle
\tableofcontents
\section{Introduction}
In this paper, we investigate the relationship between the
$L^p$ integrability of functions and the curvature geometry
of their Fourier support. Agranovsky and Narayanan
\cite{AN04} proved that if $f\in L^p(\R^n)$ and
$\operatorname{supp}\widehat f$ is carried by a $C^1$
manifold of dimension $0<d<n$, then $f\equiv0$, provided
$2\le p\le 2n/d$. They also proved that the threshold
$2n/d$ is sharp in general when $2d\ge n$. In this range,
there exist a $d$-dimensional manifold $M\subset\R^n$
and a nonzero function $f$ with
$\operatorname{supp}\widehat f\subseteq M$ and
$f\in L^p(\R^n)$ for every $p>2n/d$.

A natural question is whether the classical threshold
$2n/d$ can be improved uniformly for sufficiently regular
manifolds in the complementary range $2d<n$. Here the
normal spaces have larger dimension than the tangent
spaces. Does this force enough directions of slower
Fourier decay to improve the threshold using only
dimensional information, and what is the critical
dimensional barrier for such an improvement? We answer
these questions affirmatively in the following more
general theorem.

\begin{theorem}\label{thm:smooth-high-codimension}
Let $u\in L^p(\R^n)$, and suppose that
$\operatorname{supp}\widehat u$ is carried by a $C^n$
manifold $M$ of dimension $d<n$. Suppose $2\le k\le d+1$
is an integer satisfying
\[
n>kd-\frac{(k-1)(k-2)}{2}.
\]
Then $u\equiv0$ whenever
\[
2\le p\le
\frac{6n-(k-1)(k+2)}{3d-k+1}.
\]
\end{theorem}

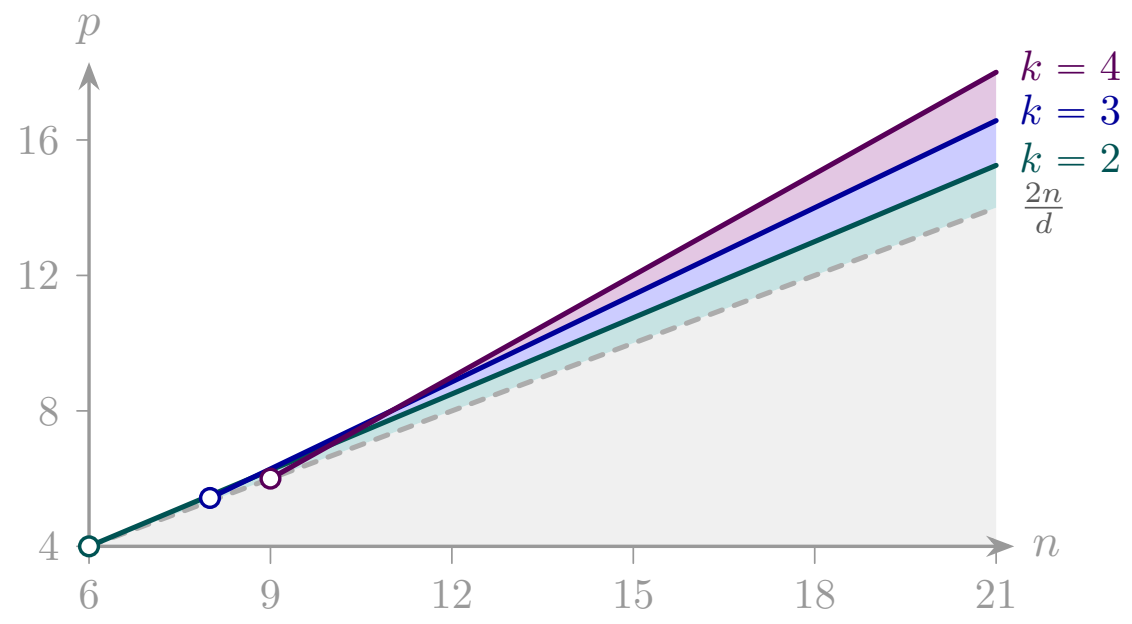
\begin{figure}[H]
\centering
\scalebox{0.8}{%
\begin{tikzpicture}[x=0.5cm,y=0.28cm,>=Stealth,line cap=round]
\def\cA{teal!65!black}\def\cB{blue!60!black}\def\cC{violet!70!black}

\fill[gray!12] (6,4) -- (21,14) -- (21,4) -- cycle;
\fill[teal!22] (6,4) -- (21,14) -- (21,{(6*21-4)/8}) -- cycle;
\fill[blue!20] ({26/3},6) -- (21,{(6*21-4)/8}) -- (21,{(6*21-10)/7}) -- cycle;
\fill[violet!22] (11,8) -- (21,{(6*21-10)/7}) -- (21,{(6*21-18)/6}) -- cycle;

\draw[->,thick,gray!80] (6,4) -- (21.3,4) node[right,font=\small] {$n$};
\draw[->,thick,gray!80] (6,4) -- (6,18.3) node[above,font=\small] {$p$};
\foreach \x in {6,9,12,15,18,21}
  \draw[gray!80] (\x,4) -- (\x,3.5) node[below,font=\footnotesize] {$\x$};
\foreach \y in {4,8,12,16}
  \draw[gray!80] (6,\y) -- (5.8,\y) node[left,font=\footnotesize] {$\y$};

\draw[gray!65,very thick,dashed] (6,4) -- (21,14);
\node[gray!70!black,font=\footnotesize,anchor=west] at (21.1,14) {$\frac{2n}{d}$};

\draw[\cA,very thick,domain=6:21,samples=2] plot (\x,{(6*\x-4)/8});
\draw[draw=\cA,fill=white,thick] (6,4) circle (2.2pt);
\node[\cA,font=\footnotesize,anchor=west] at (21.1,15.5) {$k=2$};

\draw[\cB,very thick,domain=8:21,samples=2] plot (\x,{(6*\x-10)/7});
\draw[draw=\cB,fill=white,thick] (8,{(48-10)/7}) circle (2.2pt);
\node[\cB,font=\footnotesize,anchor=west] at (21.1,16.9) {$k=3$};

\draw[\cC,very thick,domain=9:21,samples=2] plot (\x,{(6*\x-18)/6});
\draw[draw=\cC,fill=white,thick] (9,6) circle (2.2pt);
\node[\cC,font=\footnotesize,anchor=west] at (21.1,18.2) {$k=4$};
\end{tikzpicture}%
}
\caption{The exponent $\frac{6n-(k-1)(k+2)}{3d-k+1}$ of
Theorem~\ref{thm:smooth-high-codimension} for $d=3$ and $k=2,3,4$ (so $k\le d+1$),
against the Agranovsky--Narayanan exponent $2n/d$ (dashed). Larger $k$ gives a larger gain once $n$ is large enough.}
\label{fig:kfamily}
\end{figure}
\FloatBarrier

\begin{remark}
Taking $k=2$ gives
\[
\frac{6n-4}{3d-1}
=
\frac{2n}{d}
+\frac{2(n-2d)}{d(3d-1)},
\]
a strict improvement over the classical threshold for every
$C^n$ manifold whenever the codimension is strictly larger
than the dimension.
The correction term $\binom{k-1}{2}$ in the dimensional
hypothesis comes from the symmetry of the normal Hessians.
In particular, the theorem applies whenever $n>kd$,
and for $k\ge3$ also covers additional cases with $n\le kd$.
The regularity assumption can be weakened using additional
technical arguments, which we do not pursue here.
\end{remark}
The nonexistence result above also gives a necessary condition
for Fourier extension estimates. Let $M\subset\mathbb{R}^n$
be a smooth manifold of dimension $d$. Fix a nonzero
nonnegative smooth density $d\sigma$ compactly supported
in a coordinate patch of $M$. Define the extension operator by
\[
E_\sigma g(\xi)
=\int_M e^{ i\xi\cdot x}g(x)\,d\sigma(x).
\]
In its dual formulation, the Fourier restriction problem is
to determine the range of $(p,q)$ for which
\begin{equation}\label{eq:extension-estimate}
\|E_\sigma g\|_{L^q(\mathbb{R}^n)}
\lesssim_{M,\sigma,p,q}
\|g\|_{L^p(M,d\sigma)}.
\end{equation}
Testing on a nonzero smooth bump supported in a small ball
and using the classical integrability obstruction gives the
necessary condition $q>2n/d$.

Mockenhaupt conjectured that extension estimates can hold
for $q$ arbitrarily close to $2n/d$ only when the codimension
is at most half the ambient dimension, and proved this in
odd ambient dimensions \cite{Mockenhaupt96}. The following corollary proves
the conjecture in all dimensions and quantifies the gap
above $2n/d$ when $n>2d$.

\begin{corollary}\label{cor:mockenhaupt-codimension}
Suppose $n>2d$. Then, for every $1\leq p\leq\infty$, the extension
estimate \eqref{eq:extension-estimate} holds only if
\[
q>\frac{6n-4}{3d-1}
 =\frac{2n}{d}+\frac{2(n-2d)}{d(3d-1)}.
\]
\end{corollary}

\begin{proof}
Taking $g\equiv1$, we have $g\in L^p(M,d\sigma)$ for every
$1\le p\le\infty$, since $\sigma$ is finite.
As $\sigma$ is nonzero, taking $k=2$ in
Theorem~\ref{thm:smooth-high-codimension} gives
\[
E_\sigma1\in L^q(\R^n)
\implies q>\frac{6n-4}{3d-1}.
\]
This proves the necessary condition.
\end{proof}

A standard tool for obtaining necessary conditions for
Fourier restriction is the Knapp cap construction. The
corollary above shows how $p$-thinness provides an alternative
route to such conditions. An extension estimate implies
$L^q$ integrability for the Fourier transforms of smooth
densities on the manifold, while a $p$-thin result rules out
this integrability in a prescribed range. This approach is
particularly useful in higher codimension, where identifying
cap constructions adapted to the normal geometry can be
difficult. We develop this connection further in our
discussion of restriction and apply $p$-thin results to
obtain additional necessary conditions. There we also
present a more general version of the corollary,
corresponding to the full range of $k$ in
Theorem~\ref{thm:smooth-high-codimension}.
\subsection{History and Motivation}
 This problem  has precedents primarily in two areas: the study of uniqueness for convolution equations and uniqueness for homogeneous constant-coefficient PDEs. Their theorem harks back to similar results proved by Agmon--H\"ormander \cite{AH76}, Agranovsky \cite{ABCP94,ABK96}, and Thangavelu \cite{Thangavelu}. For certain degenerate hypersurfaces, Guo \cite{Guo93} obtained a stronger bound. Related problems have also been studied in other settings, including fractal \cite{SR14, DPR}, discrete \cite{ILY,BDI}, and Riemannian \cite{IMW}. Remarkably, the same dimension-based exponent $2n/d$ appears in fractal, discrete, and Riemannian analogues, with appropriate notions of dimension and integrability.

Classically, this has been known as the $p$-thin problem\cite{Guo93}, while recent literature\cite{IMW2026,BDI} has also referred to it as a spectral synthesis problem. Informally, the question can be stated as follows:

Suppose $M\subset\R^n$ is a sufficiently regular $d$-dimensional submanifold, with $0<d<n$. What is the critical exponent $p_c(M)\ge2$ such that
\[
u\in L^p(\R^n),
\qquad
\operatorname{supp}\widehat u\subseteq M,
\qquad
2\le p\le p_c(M)
\]
imply $u\equiv0$?

Theorem~\ref{thm:smooth-high-codimension} gives a lower
bound for $p_c(M)$ using only the dimension of $M$ and
the ambient dimension. For a particular manifold, however,
how much more can its geometry tell us? Consider
$M_1=\mathbb S^{n-1}$ and $M_2=\mathbb H_n$, a hyperplane
in $\R^n$. The dimension-based bound of Agranovsky and
Narayanan gives $p_c(M_i)\ge 2n/(n-1)$ in both cases.
For the sphere this is sharp. Bessel function asymptotics
show that the Fourier transform of spherical surface
measure belongs to $L^p(\R^n)$ precisely when
$p>2n/(n-1)$. For the hyperplane, on the other hand,
no nonzero function in a finite $L^p$ space can have its
Fourier transform supported on $\mathbb H_n$, thus
$p_c(\mathbb H_n)=\infty$. For the moment curve, Guo, Iosevich, Zhang, and Zorin-Kranich\cite{GIZZ23} established the critical threshold
\[
\Gamma=\{(t,t^2,\ldots,t^n):t\in[-1,1]\},
\qquad
p_c(\Gamma)=\frac{n^2+n+2}{2},
\]
which is larger than the dimension-based exponent $2n$ when $n\ge3$. Their proof uses a geometric approach involving local curvature geometry, drawing on ideas from Fourier restriction and decoupling theory. In particular, these examples bring us to a formulation of the critical exponent that records not only the dimension of the support but also its finer curvature geometry. More recently, Deodhar-Iosevich\cite{DI26} used the notion of the Fourier ratio, defined below, to introduce an additional complexity parameter, which can also encode curvature, alongside dimensional considerations. Under their hypotheses, this improves the exponent from 
\[
\frac{2n}{d}~~ \rightarrow\frac{2(n-2\kappa)}{d-2\kappa},
\]
where $\kappa$ is the Fourier-ratio decay exponent. This refinement quantifies how concentration of Fourier mass strengthens uniqueness. Their result yields sharp critical thresholds for certain classes of manifolds. These examples lead us to the following classification
question. What geometric information about $M$ is
sufficient to determine $p_c(M)$, or at least to give
a range for it? We study two classes from this point of
view. The first consists of quadratic manifolds, a
well-studied class in Fourier restriction and decoupling.
The second consists of manifolds with constant relative
nullity, for which we extend Guo's hypersurface framework
to higher codimension. We then use these results to obtain
necessary conditions for Fourier restriction and uniqueness
criteria for homogeneous constant-coefficient PDEs.

\subsection{Critical exponents for spectral synthesis}
We first consider quadratic manifolds. For a quadratic
graph, the Hessians of the graphing functions are
constant, so the eigenvalues of the normal Hessian depend
only on the normal direction. When several eigenvalues
are small, the available Fourier decay is weaker. To
estimate $p_c(S)$, we therefore need to know how large
the set of normal directions with this weaker decay can
be. This is the role of the joint sublevel estimates in
the following theorem.

We write $k_0$ and $k$ for the minimum and
maximum ranks of the normal Hessian pencil. The vector
$\mathbf a=(a_{k_0+1},\ldots,a_k)$ belongs to the cube $[0,1]^{k-k_{0}}$ and records the scales at which the
corresponding eigenvalues are small. We set
\[
K(\mathbf a)
=
k_0+\sum_{i=k_0+1}^{k}(1-a_i).
\] We state the following Theorem.
\begin{theorem}\label{Thm:QuadraticPthin}
Let $S$ be a $d$-dimensional quadratic manifold in
$\R^{n}$. Suppose that, over a compact subset,
$S$ can be reparametrized as a graph
\[
S=\{(x,\Phi(x)):x\in U,\ \Phi=(\phi_1,\ldots,\phi_{n-d})\}.
\]
Write
\[
Q(\eta)=\sum_{j=1}^{n-d}\eta_jD^{2}{\phi_j}.
\]
Suppose the minimum rank of $Q$ over
$\eta\in\R^{n-d}\setminus\{0\}$ is $k_0>0$, and its maximum
rank is $k$. Let $\Gamma_i$ denote the absolute values
of the eigenvalues of $Q$, in descending order. Suppose we have the sublevel set estimates
\[
R^{-\beta_1(\mathbf a)}
\lesssim
\left|
\left\{
\theta\in \S^{n-d-1}:
\Gamma_i(\theta)\lesssim_{S,n,d}R^{-a_i},
\quad k_0+1\le i\le k
\right\}
\right|
\lesssim
R^{-\beta_2(\mathbf a)}
\]
for sufficiently large $R$. Then, for any compactly supported distribution $T$ on $S$,
if $\widehat T\in L^p(\R^{n})$ for
\[
2\le p<
\sup_{\mathbf a}
\frac{2(n-d+K(\mathbf a)-\beta_1(\mathbf a))}
     {K(\mathbf a)},
\]
we have $T\equiv0$. If this supremum is attained, the
conclusion also holds at the critical endpoint. Moreover, if $T$ is a smooth compactly supported density
on $S$, then $\widehat T\in L^p(\R^{n})$ for every
\[
p>
\sup_{\mathbf a}
\frac{2(n-d+K(\mathbf a)-\beta_2(\mathbf a))}
     {K(\mathbf a)}.
\]
Consequently,
\[
\sup_{\mathbf a}
\frac{2(n-d+K(\mathbf a)-\beta_1(\mathbf a))}
     {K(\mathbf a)}
\le p_c(S)\le
\sup_{\mathbf a}
\frac{2(n-d+K(\mathbf a)-\beta_2(\mathbf a))}
     {K(\mathbf a)}.
\]
When we have exact sublevel set estimates for the maximizing $\mathbf a$, we obtain the sharp result.
\end{theorem}
\begin{remark}
    The proof of the result above builds on methods developed
by Mockenhaupt and Cao et al.\cite{CLMP}. In particular, small
eigenvalues of the normal Hessian weaken the available
Fourier decay bounds. Joint sublevel estimates measure the
size of the normal directions where several eigenvalues are
simultaneously small, allowing us to balance the weaker decay
against the size of the region where it occurs. The proof
uses the Fourier ratio to quantify the spectral complexity
of densities supported on these manifolds. A more detailed discussion of these sublevel estimates appears
in the corresponding section. In codimension two, we determine the exact critical exponent
from the ranks and eigenvalue vanishing orders of the normal
Hessian pencil. The resulting finite formula gives both
nonexistence at the critical endpoint and $L^p$ integrability
of nonzero smooth densities above it; see
Theorem~\ref{thm:codimtwo}.
\end{remark}
We next consider manifolds whose curvature
can vary from point to point, but whose flat tangent
directions retain a useful geometric structure. For
manifolds with constant relative nullity, these directions
locally form pieces of affine planes along which the
tangent space remains constant. This structure allows us
to obtain integrability bounds that reflect the number
of curved tangent directions.

This class includes quadratic manifolds, certain
perturbations of them, highly degenerate polynomial
surfaces, and nonpolynomial examples, as discussed later.
To state our next result, we make the meaning of manifolds
with constant relative nullity precise.
\begin{definition}[Relative nullity]
Let $M\subset\R^n$ be a smooth $d$-dimensional submanifold.
Fix $p\in M$ and choose a local parametrization
$X:U\subset\R^d\to M$ with $X(x_0)=p$. For each normal
vector $N\in(T_pM)^\perp$, define the second fundamental
matrix in the direction $N$ by
\[
H_N(p;X)
=
\left(
\left\langle
\frac{\partial^2X}{\partial x_i\partial x_j}(x_0),N
\right\rangle
\right)_{i,j=1}^{d}.
\]
The relative nullity of $M$ at $p$ is
\[
\nu(p)
=
\dim\left(
\bigcap_{N\in(T_pM)^\perp}\ker H_N(p;X)
\right).
\]
This number is independent of the choice of parametrization.
We say that $M$ has constant relative nullity $\nu$ if
$\nu(p)=\nu$ for every $p\in M$.
\end{definition}

\begin{theorem}\label{thm:PthinConstantNullity}
Let $M\subset\mathbb R^n$ be a smooth $d$-dimensional submanifold
with constant relative nullity $\nu$, and set $\rho=d-\nu$.
Suppose $\rho>0$. If $u\in L^p(\mathbb R^n)$ and
$\operatorname{supp}\widehat u\subseteq M$, then $u\equiv0$ whenever
\[
1\leq p\leq 2+\frac{2(n-d)}{\rho}.
\]
Moreover, suppose $M$ has a graph patch
\[
M_0=\{(x,\Phi(x)):x\in U\},
\qquad \Phi=(\phi_1,\ldots,\phi_{n-d}),
\]
on which
\[
\operatorname{rank}\left(
\sum_{j=1}^{n-d}\eta_jD^2\phi_j(x)
\right)\geq r>0
\qquad
\text{for every }x\in U
\text{ and }\eta\in\mathbb R^{n-d}\setminus\{0\}.
\]
Then there exists a nonzero $f\in C_c^\infty(M_0)$ such that
\[
\widehat{f\,d\sigma}\in L^p(\mathbb R^n)
\qquad\text{for every }p>2+\frac{2(n-d)}{r}.
\]
Consequently,
\begin{equation}\label{eq:range}
\boxed{
2+\frac{2(n-d)}{\rho}
\leq p_c(M)
\leq 2+\frac{2(n-d)}{r}.
}
\end{equation}
\end{theorem}
\begin{remark}
For hypersurfaces, constant relative nullity $\nu$ means
that exactly $\nu$ principal curvatures vanish at every
point. Familiar examples include spheres, circular cylinders,
and the smooth parts of circular cones. This class has been
studied by Guo and M\"uller\cite{Guo93, Muller82}. A more detailed discussion
and background can be found in Section\ref{RelativeNullity}.
\end{remark}

\subsection{Fourier Restriction for Manifolds of higher co-dimension}
We record implications of spectral synthesis results for
the ranges of Fourier restriction exponents. In higher
codimension, our results give a smaller universal outer
region for the extension exponents, with a further
improvement for quadratic manifolds.

Let $M\subset\mathbb R^n$ have dimension $d$ and
codimension $m=n-d\ge1$. Fix a nonzero nonnegative
density $\sigma$ compactly supported in a graph patch,
with smooth amplitude in the graph parameters. Write
$\mathcal T_\sigma(M)$ for the set of
$(x,y)=(1/r,1/s)\in[0,1]^2$ for which
\[
    \|E_{\sigma}g\|_{L^s(\mathbb R^n)}
    \lesssim \|g\|_{L^r(d\sigma)}.
\]
By duality, this is the restriction estimate from
$L^{s'}(\mathbb R^n)$ to $L^{r'}(d\sigma)$.

For manifolds of higher codimension, Cao, Li, Miao,
and Pang ask for the best extension region attainable
in fixed dimension and codimension
\cite[Question~7.15]{CLMP}. Let $D=D(d,n)$ be the least
sum of the degrees of $n$ distinct nonconstant monomials
in $d$ variables. Their degree and cap obstructions
give the trapezoid
\[
  \tau_0(d,n)=
  \left\{(x,y)\in[0,1]^2:
         y\le\frac{d}{D+d},\quad
         x+\frac{D}{d}y\le1
  \right\}
\]
as an outer bound for $\mathcal T_\sigma(M)$.
The cap obstruction gives the oblique boundary, while
the integrability of the Fourier transform of the density
determines the horizontal boundary. We improve the latter.

In the range
\[
2d<n\le d+\frac{d(d+1)}2,
\]
we have $D=2n-d$, so the horizontal condition above
is precisely $s\ge2n/d$. As a direct consequence of
Theorem~\ref{thm:smooth-high-codimension}, we obtain
a strict improvement in this range. More generally,
the following theorem applies throughout the dimensional
range allowed by that result and generalizes
Corollary~\ref{cor:mockenhaupt-codimension}.

\begin{theorem}\label{HighcoDimrestriction}
Let $M\subset\mathbb R^n$ be a $C^n$ manifold of
dimension $d$. Suppose $2\le k\le d+1$ is an integer
satisfying
\[
n>kd-\frac{(k-1)(k-2)}2.
\]
Then
\[
\mathcal T_\sigma(M)
\subseteq
\tau_0(d,n)\cap
\left\{(x,y)\in[0,1]^2:
y<\frac1{p_k(d,n)}\right\},
\]
where
\[
p_k(d,n)=
\frac{6n-(k-1)(k+2)}{3d-k+1}.
\]
\end{theorem}

The proof follows the same lines as that of
Corollary~\ref{cor:mockenhaupt-codimension}.
In the proof of
Theorem~\ref{thm:smooth-high-codimension}, terms beyond
quadratic order are treated as errors. For quadratic
manifolds, these terms are absent, allowing us to exploit
the normal Hessians more fully and obtain the stronger
dimensional bound below.

Let
\[
M_Q=\{(z,Q_1(z),\ldots,Q_m(z)):z\in U\},
\]
where $U\subset\mathbb R^d$ is open and the $Q_j$
are linearly independent real homogeneous quadratic
forms. Define
\[
 P_{\rm quad}(d,m):=
 \max_{\substack{0\le h<d,\ h\in\mathbb Z\\
                  m>dh-h(h-1)/2}}
 \left(2+\frac{2m-h(h+1)}{d-h}\right).
\]
We record the quadratic improvement as follows.

\begin{theorem}\label{QuadraticImprovement}
Let $M_Q$ be as above. Then
\[
 \mathcal T_\sigma(M_Q)\subseteq\tau_Q(d,m):=
 \left\{(x,y)\in[0,1]^2:
       y<\frac1{P_{\rm quad}(d,m)},\quad
       x+\frac{d+2m}{d}y\le1
 \right\}.
\]
\end{theorem}

Linear independence of the quadratic forms forces
\[
m\le\frac{d(d+1)}2,
\qquad\text{equivalently}\qquad
n\le d+\frac{d(d+1)}2.
\]
When $n>2d$, the quadratic bound gives a further
horizontal improvement over the general smooth bound.
Finer joint sublevel estimates, as discussed in the
previous section, can yield additional improvements
when more information about the normal Hessian pencil
is available.

\begin{figure}[ht]
\centering
\begin{tikzpicture}[>=Stealth,line cap=round,x=7cm,y=1cm,font=\small]
\def\cA{gray!70!black}\def\cB{teal!65!black}\def\cC{violet!70!black}
\def\H{4.2}   
\def\hA{3.3}\def\hB{2.6}\def\hC{1.9}   

\fill[blue!10]   (0,0) -- (1,0) -- ({1-\hA/\H},\hA) -- (0,\hA) -- cycle;
\fill[teal!22]   (0,0) -- (1,0) -- ({1-\hB/\H},\hB) -- (0,\hB) -- cycle;
\fill[violet!24] (0,0) -- (1,0) -- ({1-\hC/\H},\hC) -- (0,\hC) -- cycle;

\draw[->,thick,gray!80] (0,0) -- (1.12,0) node[right] {$x$};
\draw[->,thick,gray!80] (0,0) -- (0,4.8) node[above] {$y$};

\draw[very thick,\cA] (1,0) -- (0,\H);
\node[\cA,anchor=south west,font=\footnotesize] at (0.76,1.15) {$3x+11y=3$};
\draw[gray!80] (1,0.08) -- (1,-0.08);
\node[below,font=\footnotesize] at (1,-0.08) {$1$};
\draw[gray!80] (-0.012,\H) -- (0.012,\H);
\node[left,font=\footnotesize] at (-0.012,\H) {$\frac{3}{11}$};

\draw[very thick,dashed,\cA] (0,\hA) -- ({1-\hA/\H},\hA);
\draw[very thick,\cB]        (0,\hB) -- ({1-\hB/\H},\hB);
\draw[very thick,\cC]        (0,\hC) -- ({1-\hC/\H},\hC);
\node[circle,fill=\cA,inner sep=1.5pt] at ({1-\hA/\H},\hA) {};
\node[circle,fill=\cB,inner sep=1.5pt] at ({1-\hB/\H},\hB) {};
\node[circle,fill=\cC,inner sep=1.5pt] at ({1-\hC/\H},\hC) {};

\node[left,font=\footnotesize,\cA] at (-0.012,\hA) {$\frac{3}{14}$};
\node[left,font=\footnotesize,\cB] at (-0.012,\hB) {$\frac{4}{19}$};
\node[left,font=\footnotesize,\cC] at (-0.012,\hC) {$\frac{1}{5}$};

\node[\cA,anchor=west,font=\footnotesize] at ({1-\hA/\H+0.03},\hA) {degree bound, $s\ge\frac{14}{3}$};
\node[\cB,anchor=west,font=\footnotesize] at ({1-\hB/\H+0.03},\hB) {smooth theorem ($k=2$), $s>\frac{19}{4}$};
\node[\cC,anchor=west,font=\footnotesize] at ({1-\hC/\H+0.03},\hC) {quadratic theorem ($h=1$), $s>5$};

\draw[->,teal!70!black,thin,dashed] (0.06,{\hA-0.08}) -- node[right,pos=0.75,font=\scriptsize,text=teal!70!black] {improvement} (0.06,{\hC+0.08});
\end{tikzpicture}
\caption{Trapezoids with the common oblique side $3x+11y=3$ and top edges at heights $y=1/s$ for
$d=3$, $m=4$ ($n=7$, $D=11$).}
\label{fig:intro-trapezoid}
\end{figure}
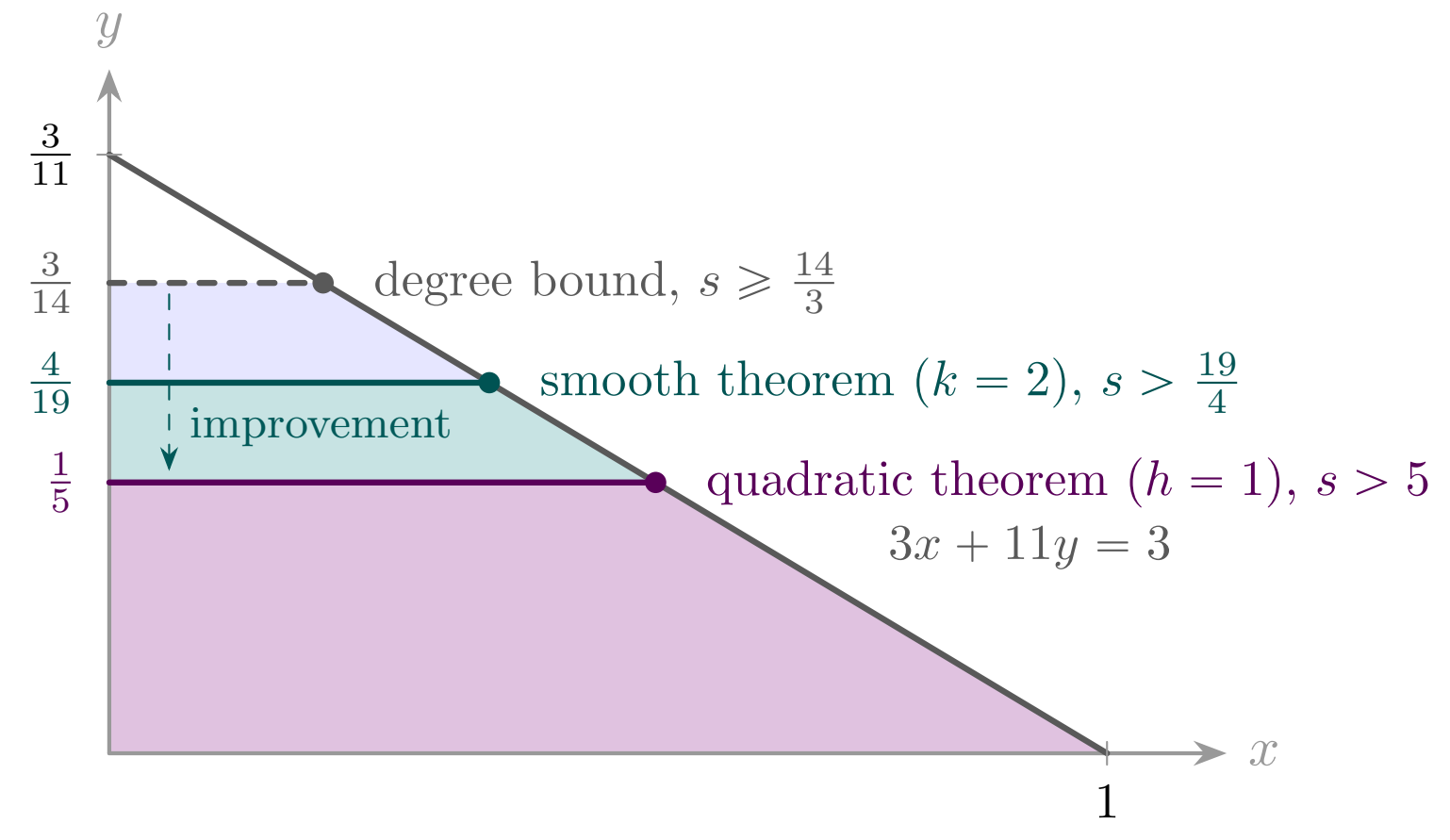

For example, take $d=3$ and $m=4$, so $n=7$ and
$D=11$. The degree bound gives $s\ge14/3$, the smooth
theorem gives $s>19/4$ ($k=2$), and the quadratic
theorem gives $s>5$ ($h=1$). Thus the successive
heights in Figure~\ref{fig:intro-trapezoid} are
$3/14$, $4/19$, and $1/5$, with the same oblique
line $3x+11y=3$.

\FloatBarrier
 We close this section by determining exactly which
dimensions admit quadratic manifolds satisfying the
Christ--Mockenhaupt condition. This classical condition was introduced by Christ
\cite{Christ82}, \cite{Christ85} and Mockenhaupt \cite{Mockenhaupt96}
as a notion of nondegeneracy in higher codimension,
in their work on Fourier restriction. More recently,
it has been used to study Fourier restriction and
decoupling for quadratic manifolds \cite{GO22,CLMP}.

Let $n>d\ge1$, and let
$Q=(Q_1,\ldots,Q_{n-d})$ be a tuple of real quadratic
forms on $\R^d$. Write
\[
H_Q(\theta)=\sum_{j=1}^{n-d}\theta_jD^2Q_j.
\]
Recall that $Q$ satisfies the Christ--Mockenhaupt
condition if
\begin{equation}\label{eq:CM}
 \int_{\S^{n-d-1}}
 |\det H_Q(\theta)|^{-\gamma}d\theta<\infty
 \qquad\text{for every }0<\gamma<\frac{n-d}{d}.
\end{equation}
For quadratic manifolds, this condition is equivalent
to $L^p$ integrability of the Gaussian-weighted extension
for every $p>2n/d$
\cite[Theorem~1.2]{CLMP}. The CM condition also plays a central role in the
conjectural description of optimal restriction ranges.
Cao, Li, Miao, and Pang conjectured that, when $n\le2d$,
the outer bound $\tau_0(d,n)$ is attained, up to
endpoints, precisely by manifolds satisfying the CM
condition \cite[Conjecture~7.20]{CLMP}.

Cao, Li, Miao, and Pang established the existence of
CM quadratic manifolds whenever $n\le2d$ and asked
whether any exist when $n>2d$
\cite[Question~7.18 and Proposition~7.19]{CLMP}.
Mockenhaupt had excluded the latter case when $d$ is
odd, and also when $d$ is even and some normal
combination $H_Q(\theta)$ is positive definite; see
\cite[Section~7.13]{CLMP}.
As another consequence of
Theorem~\ref{thm:smooth-high-codimension}, we answer the question\cite[Question~7.18 and Proposition~7.19]{CLMP}, including the
remaining case of even $d$.

\begin{corollary}\label{cor:CM-dimensions}
There exists an $(n-d)$-tuple of real quadratic forms
on $\R^d$ satisfying \eqref{eq:CM} if and only if
$n\le2d$.
\end{corollary}

\begin{proof}
If $n\le2d$, the existence follows from
\cite[Proposition~7.19]{CLMP}.

Suppose $n>2d$, and let
\[
M_Q=\{(z,Q(z)):z\in\R^d\}.
\]
Taking $k=2$ in
Theorem~\ref{thm:smooth-high-codimension} gives
\[
p_c(M_Q)\ge
\frac{6n-4}{3d-1}
=
\frac{2n}{d}
+\frac{2(n-2d)}{d(3d-1)}
>
\frac{2n}{d}.
\]
On the other hand, consider the Gaussian-weighted
extension
\[
u(x,t)=
\int_{\R^d}
e^{-i(x\cdot z+t\cdot Q(z))}
e^{-|z|^2/2}\,dz,
\qquad
(x,t)\in\R^d\times\R^{n-d}.
\]
We have $u(0,0)>0$ and
$\operatorname{supp}\widehat u\subseteq M_Q$.
If $Q$ satisfies the CM condition, then
\cite[Theorem~1.2]{CLMP} gives
\[
u\in L^p(\R^n)
\qquad\text{for every }p>\frac{2n}{d}.
\]
Choosing
\[
\frac{2n}{d}<p\le\frac{6n-4}{3d-1}
\]
contradicts
Theorem~\ref{thm:smooth-high-codimension}.
Thus no CM tuple exists when $n>2d$.
\end{proof}
CM condition as a measure of non-degeneracy and its relationship with spectral synthesis exponent $p$ has been studied in detail in forthcoming section.
\subsection{$L^p$ uniqueness for PDEs}

Uniqueness for linear homogeneous PDEs with constant
coefficients and the $L^p$ decay asymptotics of their
solutions have been studied by Agmon--H\"ormander\cite{AH76},
particularly in the case $p=2$. The investigation starts
with the following observation. Suppose $P$ is a real
multivariate polynomial, and write
\[
D=i
\left(
\frac{\partial}{\partial x_1},\ldots,
\frac{\partial}{\partial x_n}
\right).
\]
We consider
\[
P(D)u=0,
\qquad u\in\mathcal S'(\R^n).
\]
Taking the Fourier transform gives
\[
P(\xi)\widehat u=0
\implies
\operatorname{supp}\widehat u
\subseteq\{\xi\in\R^n:P(\xi)=0\}.
\]
Thus, $L^p$ decay estimates for $u$, combined with
information about its Fourier support, can yield conditions
guaranteeing uniqueness. In particular, one may combine
H\"ormander's $L^2$ lower bounds\cite{Horm72} with H\"older's
inequality to obtain an $L^p$ uniqueness estimate based
on the dimension of the zero set of $P$. As seen in the
previous discussion, a dimensional estimate may miss
finer geometric information.

For example, consider
\[
P_1(\xi)=\xi_1^2+\xi_2^2+\xi_3^2-1,
\qquad
P_2(\xi)=\xi_1^2+\xi_2^2-1.
\]
Their zero sets are the sphere and cylinder, respectively,
both of codimension one in $\R^3$. The classical
dimension-based criterion shows that if $u_i$ solves
$P_i(D)u_i=0$ and $u_i\in L^p(\R^3)$ for $1\le p\le3$,
then $u_i\equiv0$.
For $P_1$, the Fourier transform of spherical surface
measure gives a nonzero solution, and Bessel function
asymptotics show that it belongs to $L^p(\R^3)$ for every
$p>3$. For $P_2$, however, Guo's result \cite{Guo93}
shows that every nonzero solution fails to belong to
$L^p(\R^3)$ for $1\le p\le4$.

We therefore formulate an $L^p$ uniqueness criterion that
takes finer curvature geometry into account alongside
dimensional information. We also consider solutions with
a prescribed annular $L^p$ growth rate. This leads us to
an annular version of the Fourier ratio introduced in
\cite{ILPY25}, together with the decay exponent studied in\cite{DI26}.
The annular formulation is adapted to the growth
conditions appearing in H\"ormander's uniqueness theorem.

\begin{definition}
Let $\sigma_M$ be surface measure on a smooth
$d$-dimensional manifold $M\subset\R^n$, and let
$0\ne f\in C_c^\infty(M)$. Set
\[
\Lambda(R)=\{x\in\R^n:R\le |x|\le2R\},
\qquad R>0.
\]
We define the \textbf{Annular Fourier ratio} by
\[
\operatorname{FR}_{\Lambda,R}(f)
=
R^{-n/2}
\frac{\|\widehat{f\,d\sigma_M}\|_{L^1(\Lambda(R))}}
     {\|\widehat{f\,d\sigma_M}\|_{L^2(\Lambda(R))}}.
\]
Moreover, we define the \textbf{Annular decay exponent} by
\[
\kappa_\Lambda(f)
=
\sup\left\{
k\in\R:
\limsup_{R\to\infty}
R^k\operatorname{FR}_{\Lambda,R}(f)=0
\right\}.
\]
\end{definition}

To capture geometric information common to all nonzero
smooth densities on the manifold, we introduce the
following exponent.

\begin{definition}
With the notation above, the
\textbf{Fourier-ratio exponent of $M$} is defined by
\[
\kappa(M)
=
\inf\left\{
\kappa_\Lambda(f):
0\ne f\in C_c^\infty(M)
\right\}.
\]
\end{definition}

The exponent $\kappa(M)$ is positive for several classes
of manifolds. We now use this geometric quantity to state
an $L^p$ uniqueness result. The growth condition below
is satisfied with $s=0$ by every global $L^p$ solution.

\begin{theorem}\label{thm:PDEUniqueness}
Let $P(D)$ be a partial differential operator with
constant coefficients, and let
\[
\{\xi\in\R^n:P(\xi)=0\}=\bigcup_j V_j
\]
be the stratification of its real zero set, where $V_j$
has dimension $d_j$ and Fourier-ratio exponent
$\kappa_j=\kappa(V_j)$. Suppose that
$u\in\mathcal S'(\R^n)\cap L^p_{\mathrm{loc}}(\R^n)$
satisfies
\begin{equation}\label{eq:Lp-decay}
\liminf_{R\to\infty}R^{-s}
\int_{R<|x|<2R}|u(x)|^p\,dx=0
\end{equation}
for some $s<\min_j(n-2\kappa_j)$. If
\[
2\le p<
\min_j \sigma_{j},
\]
then $P(D)u=0$ implies $u\equiv0$,
where 
\[
        \sigma_j=
        \begin{cases}
            \displaystyle
            \frac{2(n-s-2\kappa_j)}
                 {d_{j}-2\kappa_j},
            &2\kappa_{j}<d_{j},\\[3ex]
            \infty,
            & \text{otherwise}.
        \end{cases}
    \]
\end{theorem}

A feature of this theorem is that it takes into account
the finer geometry of each stratum through its Fourier-ratio
exponent. Classical dimension-based criteria are governed
by the largest stratum dimension. Here, a well-curved
lower-dimensional stratum may impose a stricter restriction
on $p$ than a flat higher-dimensional stratum.

For the wave equations considered by Guo
\cite[Theorem~4]{Guo93}, taking $s=0$ along with \ref{thm:PthinConstantNullity} recovers his critical
exponent and the uniqueness range upto the endpoint.
The theorem extends this curvature-sensitive approach to
more general polynomial zero sets. We in fact prove a
more general result, whose statement and proof, together
with further discussion and examples, can be found in
the subsequent sections. 
\subsection{Outline of the rest of the paper}
In the next section, we clarify our conventions and notation and record
the definitions used in the subsequent work. We then focus on specific
higher-codimensional geometries, in particular quadratic manifolds,
non-CM manifolds, and manifolds with constant relative nullity. The goal
of this section is to understand the interplay between curvature and
the decay enforced by higher-codimensional geometry, and to compute
explicit consequences of our spectral synthesis theorems.

This is followed by a detailed discussion of the $L^p$-uniqueness
problem, including a general formulation, a corollary extending the
work of Guo \cite{Guo93}, and a closing example. We then present the
proofs of  theorems in the order in which they were stated,
with the proofs of the lemmas in the following section. We have
included an informal overview at the start of the longer proofs and
recorded some intermediate implications to help the reader follow
the arguments.
\section{Notations and definitions}\label{Backround}
\textbf{Conventions:}
Unless stated otherwise, we use the Fourier transform
convention
\[
\widehat f(\xi)
:=
\int_{\R^n}e^{ix\cdot\xi}f(x)\,dx,
\qquad f\in L^1(\R^n).
\]
For a finite Borel measure $\mu$, we similarly write
\[
\widehat\mu(\xi)
:=
\int_{\R^n}e^{ix\cdot\xi}\,d\mu(x).
\]
We extend the Fourier transform to tempered distributions
by duality. With this convention, Fourier inversion and
Plancherel's identity for Schwartz functions take the forms
\[
f(x)
=
(2\pi)^{-n}\int_{\R^n}
e^{-ix\cdot\xi}\widehat f(\xi)\,d\xi,
\qquad
\|\widehat f\|_2^2=(2\pi)^n\|f\|_2^2.
\]
Moreover,
\[
\widehat{i\partial_{x_j}f}(\xi)
=
\xi_j\widehat f(\xi).
\]

For a $C^3$ function $\Phi:\R^n\to\R$, we denote its
Hessian by
\[
D^2\Phi(x)
=
\left(
\frac{\partial^2\Phi}{\partial x_i\partial x_j}(x)
\right)_{1\le i,j\le n}.
\]

For nonnegative quantities $A$ and $B$, we write
$A\lesssim B$ if $A\le CB$ for some constant $C>0$.
Such constants may depend on fixed dimensional and
geometric data, but are independent of scale and
frequency variables unless stated otherwise.
We write $A\lesssim_{a,b}B$ to indicate dependence on
additional parameters $a$ and $b$.
We usually specify these dependencies at their first
occurrence and subsequently suppress them when they
are clear from the context. We write $A\gtrsim B$ for
the reverse inequality, and $A\approx B$ when both
$A\lesssim B$ and $B\lesssim A$ hold.

\textbf{Fourier Ratio:}
Let $\mu$ be a finite, compactly supported Borel measure
on $\R^n$, and let $0\ne f\in L^2(\mu)$. Fix
$\psi\in C_c^\infty(\R^n)$ with $\int_{\R^n}\psi=1$,
and set
\[
\psi_\delta(x)=\delta^{-n}\psi(x/\delta),
\qquad \delta>0.
\]
The family $(\psi_\delta)_{\delta>0}$ forms an
approximation to the identity.

For $R\ge1$, we define the regularized Fourier norms by
\begin{equation}\label{eq:Xdef}
X_{p,\mu,R}(f)
:=
\left(
R^{-n}\int_{\R^n}
\left|\widehat{(f\mu)*\psi_{R^{-1}}}(\xi)\right|^p
\,d\xi
\right)^{1/p},
\qquad 1\le p<\infty.
\end{equation}
The \textbf{Fourier ratio} of $f$ with respect to $\mu$
at scale $R$ is
\begin{equation}\label{eq:FRdef}
\operatorname{FR}_{\mu,R}(f)
:=
\frac{X_{1,\mu,R}(f)}{X_{2,\mu,R}(f)}
=
R^{-n/2}
\frac{\|\widehat{(f\mu)*\psi_{R^{-1}}}\|_1}
     {\|\widehat{(f\mu)*\psi_{R^{-1}}}\|_2}.
\end{equation}
Both norms are finite, and the denominator is positive
for nonzero $f\in L^2(\mu)$, so the ratio is well-defined.

\begin{definition}[Fourier-ratio decay exponent]
With the notation above, we define
\[
\kappa_\mu(f)
:=
\liminf_{R\to\infty}
\frac{-\log\operatorname{FR}_{\mu,R}(f)}{\log R}.
\]
\end{definition}

The Fourier ratio measures spectral concentration at
scale $R$. For fixed Fourier $L^2$ mass, a smaller
ratio corresponds to smaller Fourier $L^1$ mass and
hence greater concentration. A more detailed study
of the Fourier ratio and its consequences for uncertainty
principles, approximation, and spectral synthesis can
be found in \cite{ILPY25,DI26}.

\textbf{Zero set of Polynomials:}
The zero set of a polynomial is not necessarily a manifold; therefore, we decompose it before discussing its codimension or introducing an associated Fourier ratio. If
$f$ is analytic on $\mathbb{C}^{n}$, then $Z(f) = \{z \in \mathbb{C}^{n} :
f(z) = 0\}$ admits a locally finite stratification \cite[Chapter 1]{Chirka89}
$$Z(f) = \bigcup_{j} M_{j},$$
with the $M_{j}$ pairwise disjoint complex manifolds; a corresponding
stratification of the real zeros,
$$Z(f) \cap \R^{n} = \bigcup_{i} V_{i},$$
follows.
\section{Geometry of High codimensional Manifolds}
In the following subsections, we study classes of
higher-codimensional manifolds in greater detail,
exploring the consequences of the preceding theorems
and the role of curvature in determining $p_c$.
\subsection{Quadratic Manifolds}
In this subsection, we first record two consequences of
Theorem~\ref{Thm:QuadraticPthin} that give bounds for
$p_c$ through explicit joint sublevel estimates.
The first gives a crude range for $p_c$ of a quadratic
manifold using only the minimum and maximum ranks of
its normal Hessian pencil.
\begin{corollary}[Explicit rank bounds]\label{cor:rank}
Let
\[
S=\{(x,\Phi(x)):x\in U\subseteq\R^d\}
\subset\R^{n}
\]
be a quadratic graph, where
$\Phi=(\phi_1,\ldots,\phi_{n-d})$, and write
\[
Q(\eta)=\sum_{j=1}^{n-d}\eta_jD^2\phi_j.
\]
Suppose that the minimum and maximum ranks of $Q(\eta)$
over $\eta\ne0$ are $k_0\ge1$ and $k$, respectively.
Then
\begin{equation}\label{eq:rankbound}
\boxed{
2+\frac{2(n-d)}{k}
\le p_c(S)\le
2+2\left(\frac{n-d-1}{k_0}+\frac1k\right).
}
\end{equation}
In particular, every nonzero smooth compactly supported
density $a\,d\sigma_S$ satisfies
\[
\widehat{a\,d\sigma_S}\in L^p(\R^{n})
\qquad\text{for every}\qquad
p>2+2\left(\frac{n-d-1}{k_0}+\frac1k\right).
\]
If $k_0=k$, then
\[
p_c(S)=2+\frac{2(n-d)}{k},
\]
with nonexistence at the critical endpoint.
\end{corollary}

\begin{proof}
At $\mathbf a=0$, the sublevel set in
Theorem~\ref{Thm:QuadraticPthin} is the whole sphere,
provided the implicit constant is chosen sufficiently
large. Thus $\beta_1(0)=0$ and $K(0)=k$, giving
\[
p_c(S)\ge2+\frac{2(n-d)}{k}.
\]

To obtain the upper bound, we estimate the joint
sublevel sets using a nonzero $k\times k$ minor $P$
of $Q$. On the joint sublevel set corresponding to
$\mathbf a=(a_{k_0+1},\ldots,a_k)$, we have
\[
|\det P(\theta)|
\le\prod_{i=1}^{k}\Gamma_i(\theta)
\lesssim_{S,n}
R^{-\sum_{i=k_0+1}^{k}a_i}.
\]
Observe, $\det P(\theta)$ is a homogenous polynomial of degree $k$, and
\[\{\theta:\prod_{i=1}^{k}\Gamma_i(\theta)\leq \varepsilon\} \subseteq\{\theta:|\det P(\theta)|\leq \varepsilon\}. \]
Using the Polynomial sublevel estimate
\cite[Theorem~7.1]{CCW99}, followed by above set inclusion yeilds
\[
|\{\theta:\prod_{i=1}^{k}\Gamma_i(\theta)\leq R^{-\sum_{i=k_0+1}^{k}a_i}\}| \le\ R^{-\frac{\sum_{i=k_0+1}^{k}a_i}{k}}.
\]
Thus we may set
\[B_{2}(\mathbf{a})=\frac{\sum_{i=k_0+1}^{k}a_i}{k}.\]
To optimize across $\mathbf{a}$, set $s=\sum_{i=k_0+1}^{k}a_i$, so that
$K(\mathbf a)=k-s$ and $0\le s\le k-k_0$. Then
\[
\frac{n-d-\beta_2(\mathbf a)}{K(\mathbf a)}
=
\frac{n-d-s/k}{k-s}
=
\frac1k+\frac{n-d-1}{k-s}
\le
\frac1k+\frac{n-d-1}{k_0}.
\]
Equality holds at $a_{k_0+1}=\cdots=a_k=1$, where
$K(\mathbf a)=k_0$. The smooth-density conclusion of
Theorem~\ref{Thm:QuadraticPthin} now gives the upper
bound and the asserted $L^p$ membership.
\end{proof}
In the above case, the upper bound above can be improved when finer joint
sublevel estimates are available, for instance through
more precise information about the vanishing of the
$k\times k$ minors of the normal Hessian pencil.
In codimension two, the normal direction is described
by a single angular parameter, and the pencil drops
rank at only finitely many directions. The vanishing
orders of the eigenvalues near these directions give
matching lower and upper joint sublevel estimates.
This allows us to determine the sharp critical exponent
for quadratic manifolds of codimension two.
\begin{theorem}[Sharp exponent in codimension two]\label{thm:codimtwo}
Let $S$ be a $d$-dimensional quadratic manifold in
$\R^{d+2}$. Suppose that, over a compact subset,
$S$ can be reparametrized as a graph
\[
 S=\{(x,\Phi(x)):x\in U,\ \Phi=(\phi_1,\phi_2)\}.
\]
Write
\[
 H(\eta)=\eta_1H_1+\eta_2H_2.
\]
If $H_{1}$ and $H_{2}$ are linearly dependent then $p_{c}(S)= \infty$. Suppose $k$ is the maximum rank of $H(\eta)$. Set
\[
 \mathcal C
 =
 \{\eta_1=0\}
 \cup
 \left\{
 t=\frac{\eta_2}{\eta_1}:
 \eta_1\ne0,\ \operatorname{rank}H(\eta)<k
 \right\}.
\]
For $r\in\mathcal C$, let $k_r$ denote the rank of the pencil
in that direction, and set $m_r=k-k_r$. Let
$\lambda_{r,1}(t),\ldots,\lambda_{r,m_r}(t)$
denote the eigenvalue branches that vanish in that direction
without being identically zero. Then we may write
them as
\[
 \lambda_{r,j}(t)=(t-r)^{b_{r,j}}h_{r,j}(t),
 \qquad
 h_{r,j}\text{ real analytic},\quad h_{r,j}(r)\ne0,
 \qquad 1\le j\le m_r.
\]
For the direction $\eta_1=0$, interchange $H_1$ and $H_2$
and define the corresponding ranks and vanishing orders
in the same way near $t=0$. Moreover, with $v_+=\max\{v,0\}$,
\begin{equation}\label{eq:codimformula}
 \boxed{
 p_c=2+2\max\left\{
 \frac2k,\;
 \max_{\substack{r\in\mathcal C\\1\le j\le m_r}}
 \frac{\displaystyle 2-\frac1{b_{r,j}}}
 {\displaystyle k_r+
    \sum_{\ell=1}^{m_r}
    \left(1-\frac{b_{r,\ell}}{b_{r,j}}\right)_+}
 \right\}.
 }
\end{equation}
Additionally for surface measure $\sigma_S$ on $S$ and any 
nonzero $a\in C_c^\infty(U)$, we have 
\[
 \widehat{a\,d\sigma_S}\in L^p(\R^{d+2})
 \quad\Longleftrightarrow\quad p>p_c,
 \qquad 2<p<\infty.
\]
\end{theorem}
\begin{proof}
Suppose $H_{1}$ and $H_{2}$ are linearly independent. We compute the joint sublevel estimates for the vanishing eigenvalues, substitute them into Theorem~\ref{Thm:QuadraticPthin}, and then optimize the resulting exponent to obtain the sharp threshold.

Away from $\eta_1=0$, write
\[
H(\eta)
=\eta_1\left(H_1+tH_2\right),
\qquad
t=\frac{\eta_2}{\eta_1}.
\]
Since $k$ is the maximum rank of the pencil, we may choose a nonzero $k\times k$ minor $P(t)$ of $H_1+tH_2$. Every finite direction where the rank drops below $k$ is a  root of $P(t)$. There are therefore only finitely many such directions, and we focus on small neighborhoods of them to determine the additional loss of decay. Fix a  bad direction $r$, and set
\[
k_r=\operatorname{rank}(H_1+rH_2),
\qquad
m_r=k-k_r.
\]
By Rellich's theorem \cite{Rellich69}, the symmetric pencil $H_1+tH_2$ admits real-analytic eigenvalue branches near $r$. Exactly $m_r$ of these branches vanish at $r$ without being identically zero. In particular for eigen values
$\lambda_{r,1}(t),\ldots,\lambda_{r,m_r}(t)$ we may write
\[
\lambda_{r,j}(t)
=(t-r)^{b_{r,j}}h_{r,j}(t),
\qquad 1\le j\le m_r,
\]
where $b_{r,j}\ge1$ is an integer and $h_{r,j}$ is real analytic with $h_{r,j}(r)\ne0$. Thus, after shrinking the neighborhood,
\[
|\lambda_{r,j}(t)|
\approx |t-r|^{b_{r,j}},
\qquad 1\le j\le m_r.
\] 
We now obtain the joint sublevel estimates. For
$\mathbf a=(a_1,\ldots,a_{m_r})\in[0,1]^{m_r}$, let
\[
 E_r(R,\mathbf a)
 =
 \left\{
 t: |t-r|<\varepsilon_r,\quad
 |\lambda_{r,j}(t)|\le C R^{-a_j}
 \ \text{for every }j
 \right\},
\]
where $C>0$ is fixed and sufficiently large. The estimates
above imply that, up to fixed multiplicative constants,
these simultaneous inequalities are equivalent to
\[
 |t-r|
 \lesssim R^{-a_j/b_{r,j}},
 \qquad 1\le j\le m_r.
\]
Thus the smallest of these scales determines the size of
their intersection:
\[
 \min_j R^{-a_j/b_{r,j}}
 =
 R^{-\beta_r(\mathbf a)},
 \qquad
 \beta_r(\mathbf a)
 :=
 \max_{1\le j\le m_r}\frac{a_j}{b_{r,j}}.
\]
In particular, we may choose $A_r,B_r>0$, independent of
$R\ge1$ and $\mathbf a$, with $A_r<\varepsilon_r$, such that
\[
 \left\{|t-r|\le A_rR^{-\beta_r(\mathbf a)}\right\}
 \subseteq
 E_r(R,\mathbf a)
 \subseteq
 \left\{|t-r|\le B_rR^{-\beta_r(\mathbf a)}\right\}.
\]
For the lower inclusion, the inequality
$\beta_r(\mathbf a)b_{r,j}\ge a_j$ for every $j$ ensures that
a sufficiently small choice of $A_r$ makes all the defining
sublevel inequalities hold simultaneously. Consequently,
\begin{equation}\label{eq:codimtwo-sublevels}
 |E_r(R,\mathbf a)|
 \approx R^{-\beta_r(\mathbf a)}.
\end{equation}
Under the angular parametrization
$\theta(t)=(1,t)/\sqrt{1+t^2}$, we have
$d\theta=dt/(1+t^2)\sim_r dt$ near $r$, so the same estimate
holds for angular measure. Thus the lower and upper
sublevel exponents in the  theorem coincide:
\[
 \beta_1(\mathbf a)=\beta_2(\mathbf a)
 =\beta_r(\mathbf a).
\]

The $k_r$ eigenvalues bounded away from zero contribute
$k_r$ to the effective rank, so the corresponding quantity
in the  theorem is
\[
 K_r(\mathbf a)
 =
 k_r+\sum_{j=1}^{m_r}(1-a_j).
\]
We order the $b_{r,j}$ increasingly, so that the corresponding
singular values are decreasing up to fixed factors.

We repeat the argument at all finite bad directions. As for the
 the direction $\eta_1=0$ we may write
 \[H(\eta)=\eta_{2}(H_{2}+uH_{1}) \qquad \mathrm{where}~~ u=\frac{\eta_{1}}{\eta_{2}}\]
 and repeat the above argument in the neighborhood of $u=0$. Away from the finitely many bad directions, all $k$
nonzero singular values are uniformly bounded below.
That region has sublevel exponent zero and effective
rank $k$. Applying the Theorem\ref{Thm:QuadraticPthin} on these finitely
many angular regions, with $n-d=2$, gives
\begin{equation}\label{eq:codimtwo-optimization}
 \frac{p_c-2}{2}
 =
 \max\left\{
 \frac2k,\;
 \max_{r\in\mathcal C}
 \sup_{\mathbf a\in[0,1]^{m_r}}
 \frac{2-\beta_r(\mathbf a)}{K_r(\mathbf a)}
 \right\}.
\end{equation}
where supremum is over
ordered exponent vectors. It remains to evaluate this optimization. Put
\[
 s_r^*=\max_j\frac1{b_{r,j}},
 \qquad
 \kappa_r(s)=k_r+\sum_{j=1}^{m_c}(1-sb_{r,j})_+.
\]
For fixed $s=\beta_r(\mathbf a)$, the choice
$a_j=\min\{sb_{r,j},1\}$ minimizes $K_r(\mathbf a)$.
Hence
\[
 \sup_{\mathbf a\in[0,1]^{m_r}}
 \frac{2-\beta_r(\mathbf a)}{K_r(\mathbf a)}
 =
 \max_{0\le s\le s_r^*}\frac{2-s}{\kappa_r(s)}.
\]
Between consecutive breakpoints $s=1/b_{r,j}$,
$\kappa_r(s)$ is positive and affine, so the ratio is
monotone or constant. Its maximum therefore occurs at
$s=0$ or at a breakpoint. The value at $s=0$ is $2/k$;
evaluating the remaining candidates in
\eqref{eq:codimtwo-optimization} gives
\eqref{eq:codimformula}. As the lower and upper sublevel set estimates coincide, the optimized estimates above supply matching bounds for $p_c$ and thus yield the sharp exponent in \eqref{eq:codimformula}.

If $H_1$ and $H_2$ are linearly dependent, then, after
interchanging them if necessary, there exists $c\in\R$ such that
\[
\operatorname{rank}H(c,1)=0.
\]
After removing affine terms, the change of variables
$t=\eta_1-c\eta_2$ leaves the complementary normal frequency
variable free. Thus the Fourier transform of a nonzero smooth
density has no decay in that variable and cannot belong to
$L^p$ for any finite $p$, thus $p_{c}(M)= \infty$. Finally, for $0\ne a\in C_c^\infty(\R^{d+2})$, Theorem~\ref{Thm:QuadraticPthin} gives
\[
\widehat{a\,d\sigma_S}\in L^p(\R^{d+2})
\qquad\text{for every }p>p_c.
\]
\end{proof}
\begin{remark}
    The sub-level set estimates belong to a well-developed theory.
Uniform sublevel estimates under derivative hypotheses
were established by Carbery--Christ--Wright
\cite{CCW99}, while polynomial distribution inequalities
were studied by Carbery--Wright \cite{CW2001}.
For more precise exponents, Newton-polyhedron and
resolution-of-singularities methods provide sharp
estimates in broad classes of analytic examples; see
Greenblatt \cite{Greenblatt2010}.
Related questions concerning growth and stability of
real-analytic integrals were developed by
Phong--Stein--Sturm \cite{PSS1999}.

\end{remark}
The vanishing orders in Theorem~\ref{thm:codimtwo}
also give a direct criterion for the
Christ--Mockenhaupt condition. Suppose first that $k=d$.
Near a bad direction $r$, the nonvanishing eigenvalues
are bounded away from zero, so
\[
|\det H(\theta(t))|
\approx
|t-r|^{\sum_{j=1}^{m_c}b_{r,j}}.
\]
Thus the CM condition
\[
\int_{S^1}|\det H(\theta)|^{-\gamma}\,d\theta<\infty
\qquad\text{for every }0<\gamma<\frac2d
\]
holds if and only if
\[
\sum_{j=1}^{m_r}b_{r,j}\le\frac d2
\qquad\text{at every bad direction }r.
\]
At the direction $\eta_1=0$, the same criterion uses
the vanishing orders of $H_2+tH_1$ near $t=0$.
If $k<d$, the determinant vanishes identically and
the CM condition fails.

The same vanishing orders therefore determine both
the sharp critical exponent and whether the pencil
satisfies the CM condition. We study further consequences
for manifolds failing this condition in the next section.
\subsection{Non-CM Manifolds in Co-dimension 2}
The CM condition defined earlier has implications for
decoupling and Fourier restriction. In particular, it
has been used as a nondegeneracy condition for quadratic
manifolds \cite{CLMP,GO22}. A corresponding condition for
general smooth manifolds can be defined similarly,
using the normal Hessians at each point.

In the context of spectral synthesis, we may ask whether
a manifold satisfies the CM condition if and only if
\[
p_c(M)=\frac{2n}{d}.
\]
For quadratic manifolds, Cao, Li, Miao, and Pang
\cite[Theorem~1.2]{CLMP} proved that the CM condition
holds if and only if the Gaussian-weighted extension
belongs to $L^p(\R^n)$ for every $p>2n/d$.
Together with the lower bound of Agranovsky and
Narayanan, this gives $p_c(M)=2n/d$ for quadratic
CM manifolds.The following corollary gives a partial converse in
codimension two. We show that if the CM condition
fails at every point, then the critical exponent is
strictly larger than $2n/d$. Thus, in this setting,
attaining the classical threshold requires the CM
condition to hold at some point of the manifold.
\begin{corollary}[CM failure in codimension two]
\label{cor:cm-failure-codimension-two}
Let $U\subset\mathbb R^d$ be open, $d\geq1$, and let
$\Phi=(\phi_1,\phi_2)\in C^\infty(U;\mathbb R^2)$. Put
\[
 M=\{(x,\Phi(x)):x\in U\},
 \qquad
 H_x(\theta)=\theta_1D^2\phi_1(x)+\theta_2D^2\phi_2(x).
\]
We say that CM holds at $x$ if, $\det H_{x}(\theta)$ is not identically zero and,
\[
 \int_{\S^1}|\det H_x(\theta)|^{-\gamma}\,d\theta<\infty
 \qquad\text{for every }0<\gamma<\frac2d.
\]
Let $\mathcal B\subset M$ be the set of points where the
CM condition fails and $\det H_x(\theta)$ does not vanish  identically. Suppose $T \in L^{p}(\R^{d+2})$ and $\operatorname{Supp}\widehat{T} \subseteq \mathcal{B}$
then $T \equiv 0$ provided
\[2 \leq p\leq 2+\frac{12-2/m_d}{3d-1} \qquad \mathrm{where}~~~ m_{d}= \left\lfloor\frac d2\right\rfloor+1. \]
\end{corollary}
As the proof of the above corollary closely follows
that of Theorem~\ref{thm:smooth-high-codimension}, we present it
later. In that proof, we use an equivalent formulation
of the CM condition in codimension two\cite{Christ82}.
In this setting, CM failure is equivalent to the
determinant of the normal Hessian pencil having a
root of multiplicity greater than $d/2$ in some normal
direction. This produces weaker decay near that
direction, and the accumulated contribution from
CM failure at every point forces a strict improvement
over $2n/d$.
\begin{remark}
The corollary can also be formulated when the determinant
vanishes identically at every point of the supporting set.
Suppose that the maximal rank of $H_x(\theta)$ is 
 $0<k<d$ throughout this set. In that case, an argument
similar to the proof of the above corollary applies,
using only the $k$ non identically vanishing eigenvalue factors, as in the proof of Theorem\ref{Thm:QuadraticPthin}. Consequently the rank deficiency itself yields an improvement
over the above threshold, with the new range,
\[
2\le p\le 2+\frac{12}{2d+k}.
\]
\end{remark}

In the preceding cases, we obtained estimates by studying
normal Hessian pencils and the regions where their
eigenvalues become small. As mentioned earlier, when
these regions vary with the base point, this approach
may become difficult to carry out. In the next section,
we consider manifolds with constant relative nullity
and obtain bounds for $p_c$ by focusing on the common
flat directions of the phase functions. This provides
another geometric approach to spectral synthesis without
requiring the same detailed joint sublevel analysis.

\subsection{Manifolds with Constant Relative Nullity}\label{RelativeNullity}
Manifolds with constant relative nullity have been
studied in the context of spectral synthesis\cite{Muller82}, where it refers to a different problem, which is formulated
in terms of closed ideals in the Fourier algebra.
For hypersurfaces, constant relative nullity is
equivalent to having a fixed number of nonvanishing
principal curvatures. Their Fourier-analytic properties
have been studied by Guo \cite{Guo93} in the $p$-thin
problem and by \cite{Zhu25} in the study of Fourier
dimension.

As highlighted earlier, the common kernel of the
Hessians of the graphing functions has constant
dimension. This provides a persistent structure of
joint flatness as the base point varies, which we
exploit in the Fourier estimates. A local description due to Hartman makes this
structure precise. M\"uller\cite{Muller82} effectively used this description
to obtain Fourier estimates, and we state a corollary
of his theorem below. For the general statement,
we refer the reader to \cite{Muller82}.
\begin{corollary}[M\"uller {\cite[Theorem~4.1]{Muller82}}]
\label{thm:MullerTangential}
Let
\[
M=\{(x,\Psi(x)):x\in U\subset\R^d\}\subset\R^n
\]
be a smooth graph with constant relative nullity $\nu$.
Then every $f\in C_c^\infty(M)$ satisfies
\[
\int_{\R^d}
\left|\widehat{f\,d\sigma_M}(\xi,\eta)\right|\,d\xi
\lesssim_{M,f}
(1+|\eta|)^{(d-\nu)/2},
\qquad \eta\in\R^{n-d},
\]
where $d\sigma_M$ is surface measure on $M$.
\end{corollary}
To prove Theorem\ref{thm:PthinConstantNullity} we record the higher co-dimensional variant of result due to Littman\cite{Litt63}.
\begin{proposition}\label{Prop:Uniformrankdecay}
Let $\Phi=(\phi_1,\ldots,\phi_{n-d})$, where each
$\phi_j\in C^{d+3}(\mathbb R^d)$, and let $A\in C_c^\infty(\mathbb R^d)$.
Suppose that
\[
\operatorname{rank}\left(\sum_{j=1}^{n-d}\eta_jD^2\phi_j(x)\right)
\geq N
\qquad\text{for all }x\in\operatorname{supp}(A)
\text{ and }\eta\in\mathbb R^{n-d}\setminus\{0\}.
\]
Then
\[
\left|\int_{\mathbb R^d}
e^{i(\xi\cdot x+\eta\cdot\Phi(x))}A(x)\,dx\right|
\lesssim_{A,\Phi}(1+|\xi|+|\eta|)^{-N/2}.
\]
\end{proposition}
We have not found a proof of this precise formulation
in the literature, so we include one in Section\ref{Appendix}.
for completeness. The lower bound in the Theorem\ref{thm:PthinConstantNullity} is infact attained, even in higher co-dimensions. We close the section with a non-polynomial example, implying that.
\begin{example}
Fix $\nu\ge1$ and consider
\[
M_1=\{(x,y,z,e^x,e^y):x,y\in\R,\ z\in\R^\nu\}
\subset\R^{\nu+4}.
\]
This manifold has codimension $n-d=2$, relative nullity $\nu$,
and $\rho=2$. Its reduced normal Hessian satisfies
\[
\det H(x,y,\theta)=e^{x+y}\theta_1\theta_2,
\]
so $r=1$. Choosing a product density, stationary phase on
each planar curve and rapid decay in the flat variables give
$\widehat{f\,d\sigma}\in L^p$ for every $p>4$.
Together with the lower bound, this yields
\[
p_c(M_1)=4
\]
\end{example}
In the next section, we state the general $L^p$
uniqueness theorem. Combining it with the geometric
information obtained above gives explicit improvements
in the $L^p$ uniqueness range for a large class of
manifolds.

\section{$L^{p}$ Uniqueness of PDE and Fourier Ratio}
In this section, we present a general version of theorem\ref{thm:PDEUniqueness} and its
corollaries and explain in detail the role of the
Fourier ratio in obtaining an improved $L^p$ uniqueness
range. The goal of the investigation is simple: to
derive a formula that accounts for the curvature
geometry of the manifold alongside the classical
dimensional bounds. The geometry and stratification of zero sets of
real-analytic functions are briefly discussed in
Section\ref{Backround}.

Suppose $u$ solves the constant-coefficient PDE
\[
P(D)u=0.
\]
Informally, we localize the Fourier transform of $u$
so that its support lies on a single stratum rather
than the full stratified zero set. Once this is
isolated, we use the $L^p$ integrability theory
developed above to obtain a necessary condition for
the existence of nonzero $u$. We can therefore study
uniqueness using information localized on the Fourier
side. 

The key difficulty is to lift our observations for
smooth densities to general distributions, and this
is where the main work lies. In the localized setting, we may also ask what the
prescribed $L^p$ behavior tells us about the structure
of a possible nonzero solution. For example, consider
the cylinder
\[
M=\{\xi\in\R^3:\xi_1^2+\xi_2^2=1\},
\]
and choose the smooth compactly supported density
\[
\mu=d\sigma_{S^1}\otimes a(\xi_3)\,d\xi_3,
\qquad 0\ne a\in C_c^\infty(\R).
\]
As discussed above, $\widehat\mu\in L^p(\R^3)$ for
every $p>4$. Now consider the two distributions
\[
\mu_i=\frac{\partial\mu}{\partial\xi_i},
\qquad i=1,3.
\]
Both Fourier transforms solve the equation with symbol
\[
P(\xi)=\bigl(\xi_1^2+\xi_2^2-1\bigr)^2.
\]
However, stationary phase estimates give
\[
\widehat{\mu_3}\in L^p(\R^3)
\qquad\text{for every }p>4,
\]
whereas $\widehat{\mu_1}$ belongs to no finite $L^p$
space. Differentiation in $\xi_3$ is tangential to the
cylinder and preserves the smooth-density structure.
Differentiation in $\xi_1$ has a nonzero transverse
component and changes the Fourier growth.

Thus, even for distributions with the same support,
the direction of differentiation can determine whether
their Fourier transforms have the required $L^p$
behavior. To record this information, we recall the
standard notion of local transverse derivative order.

\begin{definition}[Local transverse derivative order]
Let $M\subset\R^n$ be a smooth $d$-dimensional
submanifold, and let $p\in\operatorname{supp}T$, where
$T$ is a distribution supported on $M$ near $p$.
Choose local coordinates
\[
(y,z)\in\R^d\times\R^{n-d},
\qquad p=(0,0),
\qquad M=\{z=0\}.
\]
The local structure theorem for distributions supported
on a submanifold gives a finite representation
\[
T=\sum_{|\alpha|\le N}
T_\alpha(y)\otimes\partial_z^\alpha\delta_0(z),
\]
where the $T_\alpha$ are distributions in the
tangential variables.

The \textbf{local transverse derivative order} of $T$
along $M$ at $p$, denoted by
$\operatorname{ord}^{\perp}_{M,p}(T)$, is the least
integer $m\ge0$ such that every coefficient $T_\alpha$
with $|\alpha|>m$ vanishes in a neighborhood of $y=0$.
This number is independent of the choice of local
coordinates; see \cite[Section~5]{KV96}.
\end{definition}
We state the the following Theorem.
\begin{theorem}\label{thm:PDEConstantnullity}
Let $P=P(i\,\partial/\partial x)$ be a nonzero partial differential operator
    with constant coefficients, and let
    \[
        \{\xi\in\mathbb R^n:P(\xi)=0\}=\bigcup_j V_j
    \]
    be a finite stratification of its real zero set, where $V_j$ has dimension
    $d_j$ and annular Fourier-ratio decay exponent $\kappa_j$. Suppose that
    $u\in\mathcal S'(\mathbb R^n)\cap L^p_{\mathrm{loc}}(\mathbb R^n)$, and let
    \[
        \widehat u_j=\widehat u\cdot\chi_j
    \]
    where
    $\chi_j\in C_c^\infty(\R^{n})$ with $\operatorname{Supp}(\chi_{j}) \subseteq \Omega_{j}$ where $\Omega_j$ is an open neighborhood
    satisfying $\Omega_j\cap\operatorname{supp}\widehat u\subset V_j$.
    Assume that every such nonzero localization has transverse order
    at least $m_j$ and satisfies
    \begin{equation}\label{eq:Stratum_Lp-decay}
        \liminf_{R\to\infty}R^{-s_j}
        \int_{R<|x|<2R}|u_j(x)|^p\,dx=0
    \end{equation}
    for some $s_j< n-2\kappa_j$. If
    \[
        2\le p<\min_j\sigma_j,
    \]
    then $Pu=0$ implies $u=0$, where
    \[
        \sigma_j=
        \begin{cases}
            \displaystyle
            \frac{2(n-s_j-2\kappa_j)}
                 {d_j-2m_j-2\kappa_j},
            & m_j<\dfrac{d_j-2\kappa_j}{2},\\[3ex]
            \infty,
            & \text{otherwise}.
        \end{cases}
    \]
\end{theorem}
When we make no assumptions on local transversal orders, and assume global $L^{p}$ decay we recover Theorem\ref{thm:PDEUniqueness}. More generally, the
theorem restricts the possible structure of a nonzero
solution: the geometry of each stratum and the prescribed
localized $L^p$ behavior constrain the transverse order
of its Fourier transform. Higher-order transverse terms
must vanish wherever the supporting geometry provides
insufficient decay. The PDE uniqueness theorems, together with the annular
Fourier-ratio estimates, give uniqueness ranges for a
large class of characteristic sets. We revisit two
geometries to illustrate these consequences. The
following corollary generalizes the uniqueness part
of Guo's theorem \cite[Theorem~4]{Guo93}.

\begin{corollary}
Let $P(D)$ be a nonzero partial differential operator
with constant coefficients, and let
\[
\{\xi\in\R^n:P(\xi)=0\}=\bigcup_j V_j
\]
be a finite stratification of its real zero set.
Write $d_j=\dim V_j$, and suppose every
positive-dimensional stratum has constant relative
nullity $d_j-k_j<d_j$.
If $u\in L^p(\R^n)$ for some finite $p$ satisfying
\[
2\le p\le
2+\min_{j:\,d_j>0}\frac{2(n-d_j)}{k_j},
\]
then $P(D)u=0$ implies $u\equiv0$.
The minimum is interpreted as $\infty$ if there are
no positive-dimensional strata.
\end{corollary}

The proof follows the argument of
Theorem~\ref{thm:PDEConstantnullity}, using the tail
lower bound \ref{eq:NullityTail} to include the
critical endpoint. When smooth densities provide matching nonzero solutions
for every $p$ above the stated threshold, the uniqueness
range is sharp. In these cases, we extend Guo's full
$L^p$ uniqueness result to the corresponding class of
equations.
\begin{remark}
The scope of the theorem extends beyond the annular
Fourier ratio. As the proof Theorem\ref{thm:PDEConstantnullity} and Lemma\ref{lem:annular-testing} show,
a suitable Fourier-tail lower bound on a region
$\Omega_R\subseteq\Gamma_R$ suffices to obtain a
uniqueness range on each stratum. For quadratic
manifolds, finer joint sublevel estimates may improve
these bounds by identifying small regions of normal
directions with weaker decay. Such regions may control
the $L^p$ obstruction without dominating the annular
$L^1$ norm. The possible gain therefore depends on both
their volume and the Fourier mass they capture.
\end{remark}
We end the section with an example.
\begin{example}
In $\R^5$, consider the symbol
\[
P(\xi)=\xi_1
\underbrace{\Big[(\xi_1-1)^2
 +(\xi_4-2-\xi_2-2\xi_3)^2
 +(\xi_5-\xi_2^2)^2\Big]}_{\text{quadratic cylinder}}
\underbrace{\Big[(\xi_4-A)^2+(\xi_5-B)^2\Big]}_
           {\text{polynomial graph}}.
\]
where
\[
\begin{aligned}
A&=\xi_1^4+\xi_1^2\xi_2+2\xi_1\xi_3,\\
B&=\xi_1^6+3\xi_1^4\xi_2+6\xi_1^2\xi_2
   +4\xi_1^3\xi_3+12\xi_1\xi_3,
\end{aligned}
\]
Suppose $P(D)u=0$ and $u$ satisfies the annular growth
condition \eqref{eq:Lp-decay} with exponent $s<1$.
We determine the uniqueness range supplied by
Theorem~\ref{thm:PDEUniqueness}.
\end{example}

\begin{proof}[Solution]
The real zero set of $P$ is the union of a hyperplane,
the graph of $(A,B)$, and a quadratic cylinder. The
cylinder is disjoint from the other two components.
The hyperplane and the graph intersect along the
two-plane $\{\xi_1=\xi_4=\xi_5=0\}$, which we take as
a separate stratum. Using the annular Fourier-ratio
estimates for these geometries, we obtain
\[
\begin{array}{lcccc}
\hline
\text{Stratum} & d_j & 5-d_j & \kappa_j
  & \text{Uniqueness range}\\
\hline
\text{Hyperplane, off the intersection}
  & 4 & 1 & 2   & 2\le p<\infty\\
\text{Polynomial graph, off the intersection}
  & 3 & 2 & 1/2 & 2\le p<4-s\\
\text{Quadratic cylinder}
  & 2 & 3 & 1/2 & 2\le p<2(4-s)\\
\text{Plane--graph intersection}
  & 2 & 3 & 1   & 2\le p<\infty\\
\hline
\end{array}
\]
Taking the minimum over the strata, we conclude that
$u\equiv0$ whenever
\[
2\le p<4-s.
\]
By comparison, the largest stratum has dimension four,
so H\"ormander's $L^2$ estimate followed by H\"older's
inequality gives only
\[
2\le p\le\frac{5-s}{2}.
\]
Since $4-s>(5-s)/2$ for $s<1$, the improvement is
strict. Here it is the middle-dimensional stratum,
the graph of $(A,B)$, that determines the range of $p$.
\end{proof}
In the rest of the paper, we present proofs of result presented in the paper.
\section{Proofs of Theorems}
\subsection{Proof of Theorem\ref{thm:smooth-high-codimension}}
\begin{proof}
Suppose $u\in L^p(\R^n)$ and
$\operatorname{supp}\widehat u\subseteq M$. We argue by
induction on $k$. Set
\[
P_j=\frac{6n-(j-1)(j+2)}{3d-j+1}.
\]
At each step, we first establish the result for the range $2\leq p\le2k$. For $k=2$, the dimensional hypothesis gives $n>2d$.
The Agranovsky--Narayanan theorem therefore implies
$u\equiv0$ whenever $2\le p\le4$. Now let $k\ge3$ and assume that the theorem holds for
$k-1$. Since
\[
n>kd-\frac{(k-1)(k-2)}2
>
(k-1)d-\frac{(k-2)(k-3)}2,
\]
the induction hypothesis gives $u\equiv0$ for
$2\le p\le P_{k-1}$. This leads to two possible cases: If $2k \leq P_{k-1}$ then we have established the result in the desired range. If $2k>P_{k-1}$, then following the observation,
\[
P_k-P_{k-1}
=\frac{P_{k-1}-2k}{3d-k+1}
\]
we conclude $P_{k} \leq P_{k-1}$, and that resolves the problem for the entire range. Thus in particular it suffices to prove the result for the range $2k < p \leq P_{k}$. Suppose that $u\in L^p(\mathbb R^n)$ is nonzero, $2k < p \leq P_{k}$ and
$\operatorname{supp}\widehat u\subseteq M$. To proceed, we first explain the idea of the proof. We localise the
distribution and reduce the problem to estimating the $L^p$
tail of its Fourier transform. We then decompose a test function
along the manifold into small rectangles, chosen so that the
oscillations of order three and above are controlled. This allows
us to estimate each piece using the quadratic part of the phase.
Next, we extend the localised functions in the normal directions
and combine their estimates. Finally, we optimise over the
possible extensions, giving more weight to normal directions
where the quadratic oscillations are weaker. The abundance of
normal directions in high codimension makes this choice
effective, and the resulting bound, together with Fourier
duality and the vanishing $L^p$ tail, finishes the proof.

\medskip
\noindent\textbf{Step I: Reduction and regularization}\\
Suppose we locally write $M$ as the graph
\[
M_0=\{(x,\Phi(x))~|~x\in U\subseteq\R^d\},
\]
where $U$ is open and $\Phi:U\to\R^{n-d}$ is $C^{n}$ function.
Choose $\kappa\in C_c^\infty(\R^n)$ supported in this coordinate
neighbourhood so that
\[
T:=\kappa\widehat{u}\neq0.
\]
Then $\operatorname{supp}T\subseteq M_0$ and
$\check{T}\in L^p(\R^n)$ by Young's inequality. Let $r>n/2$ be an integer and let $\varphi\in C_c^\infty$
be supported in a fixed compact set $K$. On the graph chart, let
$S:\R^{n}\to R^{n}$ 
\[
S(x,y)=(x,y-\Phi(x)).
\]
Since $\check T\in L^p$, $p\ge2$, $\Phi\in C^n$  standard Fourier estimates
imply that $T\circ S^{-1}$ extends continuously to $C_c^n$. Thus we define the flattened distribution,
\[
\langle T\circ S^{-1},\psi\rangle
:=\langle T,\psi\circ S\rangle.
\]
for all test functions $\psi \in C_{c}^{\infty}(\R^{n})$. As $T \circ S^{-1}$ is supported on the manifold $\{y=0\}$, by the Schwartz structure theorem for distributions\cite{Schwartz50},
\[
T\circ S^{-1}
=\sum_{|\alpha|\leq m}
\nu_\alpha(x)\otimes\partial_y^\alpha\delta_0(y),
\qquad \nu_\alpha\in\mathcal{S}'(U).
\]
Choose $\alpha_0$ of maximal total order with
$\nu_{\alpha_0}\neq0$. Choose
$\gamma\in C_c^\infty(U\times\R^{n-d})$ with  $\gamma \equiv1$ on a
neighbourhood of $\operatorname{supp}T$, and define
\[
T_0:=\gamma(x,y)(y-\Phi(x))^{\alpha_0}T.
\]
By maximality of $\alpha_0$, only the term indexed by
$\alpha_0$ survives; this simplifies the action to
\[
\langle T_0,\psi\rangle
=\langle\nu,\psi(x,\Phi(x))\rangle,
\qquad
\nu=(-1)^{|\alpha_0|}\alpha_0!\,\nu_{\alpha_0}\neq0,
\]
for every $g\in C_c^\infty(\R^n)$.
Young's inequality gives $\widehat{T}_0\in L^p(\R^n)$, and we now
estimate its $L^p$ tail decay. Since $T_0\neq0$, we may choose
$a\in C_c^\infty(\R^d)$ such that $\langle\nu,a\rangle\neq0$.

Choose a nonnegative, even function
$\rho\in C_c^\infty(B(0,1))$ on $\R^d$ with $\int\rho=1$.
For $R>0$, set
\[
\rho_R(x)=R^{d}\rho(Rx),
\qquad
\nu_R=\nu*\rho_R,
\qquad
a_R=a*\rho_R.
\]
Since $\rho_R$ is an approximation to the identity
and $\rho$ is even,
\[
\langle\nu_R,a_R\rangle
\longrightarrow\langle\nu,a\rangle
\qquad\text{as }R\to \infty.
\]
We next split $a_R$ into smaller pieces so that the
higher-order Taylor remainder of the phase is uniformly bounded
on each piece.

\medskip
\noindent\textbf{Step II: Localization and the quadratic estimate}\\
Choose a smooth partition of unity
$(\chi_j)_{j\in\mathcal{J}(\delta)}$ on a fixed bounded
neighbourhood of $\operatorname{supp}a_R$, subordinate to
fixed enlargements of the cubes in the grid $\delta\mathbb Z^d$,
with $0\leq\chi_j\leq1$ and
\[
\|\partial^\beta\chi_j\|_\infty
\lesssim_\beta\delta^{-|\beta|},
\qquad \delta=R^{-1/3}.
\]
Define
\[
b_j=(a_R\chi_j)*\rho_R.
\]
Since $\rho_R$ is supported in $B(0,R^{-1})$
and $R^{-1}\leq\delta$, the supports of $b_j$ lie in
further fixed enlargements of these grid cubes, whose overlap
is bounded solely in terms of $d$, and 
\[
|\mathcal{J}(\delta)|\lesssim_{a,d}\delta^{-d}.
\] We now estimate the $L^2$ norm of $b_j$ and the $L^1$ norm
of the associated oscillatory Fourier transform. These bounds
will later be used to estimate the $L^p$ tail of $\widehat{T}_0$.
By Young's convolution inequality and the volume of the support,
\begin{equation}\label{eq:tangent-partition-square}
\|b_j\|_2^2
\leq\|a_R\chi_j\|_2^2\|\rho_R\|_1^2
\lesssim_{a,d}\delta^d.
\end{equation}
For the $L^1$ estimate, we introduce a Gaussian factor and split
the oscillatory function into two factors, one containing the
quadratic phase and the other containing the amplitude and the
higher-order phase. We use the cube size to control the
oscillation of the remainder and estimate the two factors
separately. Let $x_j$ be the centre of the $j$th cube, and write
\[
H_j(\omega)=\sum_{\ell=1}^c\omega_\ell D^2\Phi_\ell(x_j),
\qquad \omega\in\R^c,
\]
where $c=n-d$. Denote its singular values, in decreasing order,
by $\Gamma_{j,1}(\omega),\ldots,\Gamma_{j,d}(\omega)$.

For $R<|\eta|<2R$, we use the Taylor expansion
\[
\Phi(x_j+\delta z)
=\Phi(x_j)+\delta D\Phi(x_j)z
+\frac{\delta^2}{2}D^2\Phi(x_j)[z,z]
+\mathcal R_j(z).
\]
Choose a nonnegative function $A\in C_c^\infty(\R^d)$ equal
to $1$ on a fixed ball containing the supports of all the
functions $z\mapsto b_j(x_j+\delta z)$. Set
\[
\begin{aligned}
G_{j,\eta}(z)
&=\exp\left(
-\frac{|z|^2}{2}
+\frac{i\delta^2}{2}\langle H_j(\eta)z,z\rangle
\right),\\
B_{j,\eta}(z)
&=A(z)b_j(x_j+\delta z)
\exp\left(
\frac{|z|^2}{2}+i\eta\cdot\mathcal R_j(z)
\right).
\end{aligned}
\]
Thus we write
\[
b_j(x_j+\delta z)e^{i\eta\cdot\Phi(x_j+\delta z)}
=e^{i\eta\cdot\Phi(x_j)}
e^{i\delta\,\eta\cdot D\Phi(x_j)z}
G_{j,\eta}(z)B_{j,\eta}(z).
\]
Since translation, modulation, and rescaling preserve the Fourier
$L^1$ norm, Young's convolution inequality gives
\[
\bigl\|\widehat{b_j e^{i\eta\cdot\Phi}}\bigr\|_1
\lesssim_d
\|\widehat{G}_{j,\eta}\|_1
\|\widehat{B}_{j,\eta}\|_1.
\]
To estimate the $L^1$ norm of $\widehat{G}_{j,\eta}$, we follow
the Gaussian computation in \cite{CLMP}
to obtain
\[
\|\widehat{G}_{j,\eta}\|_1
=(2\pi)^{d}\det\!\left(I+\delta^4H_j(\eta)^2\right)^{1/4}
\approx_{d}\prod_{i=1}^d
\left(1+\delta^4\Gamma_{j,i}(\eta)^2\right)^{1/4}.
\]
Since $H_j(\eta)=R H_j(\eta/R)$, this gives
\[
\|\widehat{G}_{j,\eta}\|_1
\approx_d
\prod_{i=1}^d
\left(1+R\delta^2\Gamma_{j,i}(\eta/R)\right)^{1/2}.
\]
A more detailed computation can be found in proof of Theorem\ref{Thm:QuadraticPthin}.

To estimate the remaining factor, Young's convolution
inequality gives
\[
\|\partial^\beta b_j\|_\infty
\leq\|\partial^\beta(a_R\chi_j)\|_\infty
\|\rho_R\|_1
\lesssim_{a,\rho,\beta,d}\delta^{-|\beta|},
\]
where the last bound follows from the uniform derivative
bounds for $a_R$ and $\chi_j$. Consequently,
\[
\bigl\|\partial_z^\beta
[b_j(x_j+\delta z)]\bigr\|_\infty
=\delta^{|\beta|}\|\partial^\beta b_j\|_\infty
\lesssim_{a,\rho,\beta,d}1.
\]
Moreover, for every fixed multi-index $\beta$, on the fixed
rescaled support,
\begin{equation}
\label{eq:scaled-cubic-error}
\bigl|\partial_z^\beta
(\eta\cdot\mathcal R_j(z))\bigr|
\lesssim_{\Phi,\beta,d,c}R\delta^3=1.
\end{equation}
These
derivative bounds followed by Cauchy-Schwarz and Plancheral yield, 
\[
\begin{aligned}
\|\widehat{B}_{j,\eta}\|_1
&\leq
\left(\int_{\R^d}(1+|\xi|^2)^{-m}\,d\xi\right)^{1/2}
\left(\int_{\R^d}(1+|\xi|^2)^m
|\widehat{B}_{j,\eta}(\xi)|^2\,d\xi\right)^{1/2}\\
&\lesssim_d
\sum_{|\beta|\leq\lfloor d/2\rfloor+1}
\|\partial_z^\beta B_{j,\eta}\|_2
\lesssim_{a,\rho,\Phi,d,c}1,
\end{aligned}
\]
where $m=\lfloor d/2\rfloor+1>d/2$. 
Combining the two bounds,
\begin{equation}
\label{eq:localized-smooth-fourier-algebra}
\left\|\widehat{b_j e^{i\eta\cdot\Phi}}\right\|_1
\lesssim_{a,\rho,\Phi,d,c}
\prod_{i=1}^d
\left(1+R\delta^2\Gamma_{j,i}(\eta/R)\right)^{1/2}.
\end{equation}

\medskip
\noindent\textbf{Step III: Extensions of the localized tests}\\
The localised functions are defined in the tangential variables,
over $\R^d$. In particular, we may write
$a(x)=\psi(x,\Phi(x))$ for a smooth extension $\psi$ in the
normal directions near $M$. We now give a general construction
of such extensions for the $b_j$ and obtain $L^{p'}$ bounds for the
Fourier transform of the resulting extension. Put
\[
\mathcal A=\{\omega\in\mathbb R^c:1<|\omega|<2\}.
\]
Let $h_j\in C_c^\infty(R\mathcal A)$ satisfy $\int h_j=1$
for every $j$, and write $h=(h_j)_j$. Define
\[
\psi_{h,R}(x,y)
=(2\pi)^{n-d}\sum_j b_j(x)\check h_j(y-\Phi(x)).
\]
We note that
\[
\psi_{h,R}(x,\Phi(x))=(2\pi)^{n-d}\sum_j b_j(x),
\qquad\text{and}\qquad
\langle T_0,\psi_{h,R}\rangle
=\langle\nu_R,a_R\rangle.
\]
Fourier duality followed by Hölder's inequality gives
\[
|\langle T_0,\psi_{h,R}\rangle|
\leq
\|\widehat{T}_0\|_p
\|\widehat{\psi_{h,R}}\|_{p'}.
\]
Thus, to estimate the action of $T_0$, we seek an $L^{p'}$
bound for the Fourier transform of the extension $\psi_{h,R}$.
Taking the Fourier transform first in $y$ and then in $x$
\[
\widehat{\psi_{h,R}}(\xi,\eta)
=(2\pi)^{n-d}\sum_j h_j(\eta)
\widehat{b_j e^{i\eta\cdot\Phi}}(\xi),
\]
Using the bounds from
\eqref{eq:localized-smooth-fourier-algebra} to estimate $L^{1}$norm gives,
\begin{equation}
\label{eq:dual-operator-L1}
    \|\widehat{\psi_{h,R}}\|_1
    \lesssim_{a,\rho,\Phi,d,c}
    \sum_j\int_{R\mathcal A}
       h_{j}(\eta)\prod_{i=1}^d
\left(1+R\delta^2\Gamma_{j,i}(\eta/R)\right)^{1/2}\,d\eta.
\end{equation}
Let $R\delta^2=\tau^{-1}$, and set
\[
A_j(\eta)
=\prod_{i=1}^d
\left(1+\tau^{-1}\Gamma_{j,i}(\eta/R)\right)^{1/2}.
\]
The reciprocal of the above quantity is comparable to the
determinant factor appearing in the Gaussian Fourier decay
computation for quadratic manifolds\cite{CLMP,Mockenhaupt96}. Our use of this quantity
to study $L^p$ integrability is motivated by \cite{CLMP,Mockenhaupt96}, as in our treatment of the quadratic case.
For the $L^2$ estimate, Plancherel gives
\[
\|\widehat{\psi_{h,R}}\|_2^2
=(2\pi)^{2(n-d)}\int\!\!\int
\left|\sum_j b_j(x)\check h_j(y-\Phi(x))\right|^2
\,dy\,dx.
\]
Since the supports of the $b_j$ have bounded overlap,
Cauchy--Schwarz gives
\begin{align*}
\|\widehat{\psi_{h,R}}\|_2^2
&\lesssim_d
\sum_j\int |b_j(x)|^2
\left(\int|\check h_j(y-\Phi(x))|^2\,dy\right)dx\\
&=\sum_j\|b_j\|_2^2\|\check h_j\|_2^2
=(2\pi)^{n}\sum_j\|b_j\|_2^2\|h_j\|_2^2.
\end{align*}
Using \eqref{eq:tangent-partition-square}, we obtain
\begin{equation}
\label{eq:dual-operator-L2}
\|\widehat{\psi_{h,R}}\|_2
\lesssim_{a,\rho,d}
\delta^{d/2}
\left(\sum_j\int_{R\mathcal A}|h_j(\eta)|^2\,d\eta\right)^{1/2}.
\end{equation}
Weighted interpolation of the linear map
$h\mapsto\widehat{\psi_{h,R}}$ between its $L^1$ and $L^2$
bounds yields
\begin{equation}
\label{eq:weighted-dual-interpolation}
\|\widehat{\psi_{h,R}}\|_{p'}
\lesssim_{a,\rho,\Phi,d,c,p}
\delta^{d/p}
\left(
\sum_j\int_{R\mathcal A}
|h_j(\eta)|^{p'}A_j(\eta)^{2-p'}\,d\eta
\right)^{1/p'}.
\end{equation}
In the next step, we construct $h_{j}$ which satisfy the initial assumptions and minimize the integral on the right.

\medskip
\noindent\textbf{Step IV: Selecting minimizers and concluding estimates}\\
We now choose each $h_j\geq0$, with $\int h_j=1$, to minimise
\[
\int_{R<|\eta|<2R}
h_j(\eta)^{p'}A_j(\eta)^{2-p'}\,d\eta.
\]
Hölder's inequality gives
\[
1=\int_{R<|\eta|<2R}h_j(\eta)\,d\eta
\leq
\left(\int_{R<|\eta|<2R}h_j(\eta)^{p'}A_j(\eta)^{2-p'}\,d\eta\right)^{1/p'}
\left(\int_{R<|\eta|<2R}A_j(\eta)^{-(p-2)}\,d\eta\right)^{1/p}.
\]
Equality holds when the two powered factors are proportional,
which gives
\[
h_j(\eta)\propto A_j(\eta)^{-(p-2)}.
\]
Normalising the proportionality relation gives a minimiser
\[
h_j^*(\eta)
=
\frac{1_{\{R\leq|\eta|\leq2R\}}A_j(\eta)^{-(p-2)}}
{\displaystyle
\int_{R<|\zeta|<2R}A_j(\zeta)^{-(p-2)}\,d\zeta}.
\]
As $h_{j}^{*}$ is not smooth, to find suitable $h_j$, we take smooth approximations of $h_j^*$.
Choose $0\leq\theta_\ell\leq1$ in $C_c^\infty(R\mathcal A)$
with $\theta_\ell\to1$ pointwise on $R\mathcal A$, and a
nonnegative standard mollifier $\kappa_\varepsilon$. Set
\[
h_{j,\ell}
=
\frac{(\theta_\ell h_j^*)*\kappa_{\varepsilon_\ell}}
{\displaystyle\int\theta_\ell h_j^*}.
\]
For sufficiently large $\ell$, choose $\varepsilon_\ell>0$
small enough that $h_{j,\ell}$ is supported compactly inside
$R\mathcal A$ and
\[
\begin{aligned}
\int_{R<|\eta|<2R}
h_{j,\ell}(\eta)^{p'}A_j(\eta)^{2-p'}\,d\eta
&\leq
2\int_{R<|\eta|<2R}
h_j^*(\eta)^{p'}A_j(\eta)^{2-p'}\,d\eta\\
&=
2\left(
\int_{R<|\eta|<2R}
A_j(\eta)^{-(p-2)}\,d\eta
\right)^{-p'/p}.
\end{aligned}
\]
Fix such an $\ell$ and set $h_j:=h_{j,\ell}$. Then
$h_j\in C_c^\infty(R\mathcal A)$, $h_j\geq0$, and
$\int h_j=1$, so $h_j$ satisfies all the required properties. We may write,
\[
A_j(\eta)\approx_d
\tau^{-d/2}
\left(\prod_{i=1}^d
\max\{\tau,\Gamma_{j,i}(\omega)\}\right)^{1/2}.
\]
As $n-d>(k-1)d-\frac{(k-1)(k-2)}2$
we apply Lemma~\ref{lem:determinant-integral} with $k-1$
in place of $k$ and $s=(p-2)/2\ge k-1$ to obtain
\[
\int_{\mathcal A}(\prod_{i=1}^d
\max\{\tau,\Gamma_{j,i}(\omega)\})^{-(p-2)/2}\,d\omega
\gtrsim
\tau^{-(k-1)(p-2)/2+k(k-1)/2}.
\]
This gives a lower bound for the integral of
$A_j(\eta)^{-(p-2)}$. Consequently,
\[
\left(
\int_{R<|\eta|<2R}
A_j(\eta)^{-(p-2)}\,d\eta
\right)^{-p'/p}
\lesssim
R^{-cp'/p}
\tau^{-\frac{p'((d-k+1)(p-2)+k(k-1))}{2p}}.
\]
Since the number of rectangular $\delta$-boxes satisfies
$|\mathcal{J}(\delta)|\lesssim_{a,d}\delta^{-d}$,
we apply the above estimate on each localised box.
Substituting into \eqref{eq:weighted-dual-interpolation} gives
\[
\begin{aligned}
\|\widehat{\psi_{h,R}}\|_{p'}
&\lesssim
\delta^{d/p-d/p'}R^{-c/p}
\tau^{-\frac{(d-k+1)(p-2)+k(k-1)}{2p}}\\
&=
R^{\frac{3d-k+1}{6}
-\frac{6n-(k-1)(k+2)}{6p}},
\end{aligned}
\]
where we used $\delta=\tau=R^{-1/3}$.
The exponent is nonpositive precisely when
\[
p\leq\frac{6n-(k-1)(k+2)}{3d-k+1}.
\]
Thus $\|\widehat{\psi_{h,R}}\|_{p'}\lesssim1$.
Since $\widehat{\psi_{h,R}}$ is supported in
$\{(\xi,\eta):|\eta|>R\}$, Fourier duality and Hölder's
inequality give
\[
|\langle\nu_R,a_R\rangle|
=
|\langle T_0,\psi_{h,R}\rangle|
\lesssim
\|\widehat{T}_0\|_{L^p(\{|\eta|>R\})}
\|\widehat{\psi_{h,R}}\|_{p'}
\longrightarrow0.
\]
But $\langle\nu_R,a_R\rangle\to\langle\nu,a\rangle\neq0$,
a contradiction.
\end{proof}

\subsection{Proof Of Theorem\ref{Thm:QuadraticPthin}}
\begin{proof}
We divide the proof into four steps. In the first step, aided
by a direct computation for Gaussian densities, we estimate
the Fourier ratio of compactly supported smooth densities
to obtain a lower bound on the $L^p$ norm of their Fourier
transforms. In the second step, we combine this lower bound
with sublevel-set estimates to obtain a necessary condition
for $L^p$ integrability for all nonzero compactly supported
smooth densities. In Step~III, we estimate the $L^p$ norm
at each dyadic scale and sum over the scales to obtain a
sufficient condition for $L^p$ integrability of the Fourier
transforms of smooth densities. In the final step, we extend
the necessary condition from Step~II to all nonzero compactly
supported tempered distributions, adapting the strategy
of Lemma~0.9.

\textbf{Step I: Fourier Ratio for Quadratic Sub-manifolds}\\
We will compute the Fourier ratio to obtain a lower bound on the
$L^p$ norm of the Fourier transform of the surface measure.
In order to record the $L^1$ growth, we first recall the direct
Gaussian calculation in
\cite[Proposition~3.6 and Section~7.1]{CLMP}.
We briefly present their argument; the interested reader may refer
to \cite{CLMP}   for further details.
\begin{align}\label{GaussianSplit}
\left|\widehat{(e^{-|x|^2/2}\,d\sigma)}\right|
&= \left| \int_{\mathbb{R}^d}
e^{i(\xi\cdot x+x^T\overline{Q}(\eta)x/2)}
e^{-|x|^2/2}\,dx \right| \nonumber\\
&= \left| \int_{\mathbb{R}^d}
e^{i\xi\cdot x}
e^{-x^T(I-i\overline{Q}(\eta))x/2}\,dx \right| \nonumber\\
&= \left| (2\pi)^{d/2}
\det(I-i\overline{Q}(\eta))^{-1/2}
e^{-\xi^T(I-i\overline{Q}(\eta))^{-1}\xi/2}
\right| \nonumber\\
&\sim \left|\det(I-i\overline{Q}(\eta))\right|^{-1/2}
\left|e^{-\xi^T(I-i\overline{Q}(\eta))^{-1}\xi/2}\right|,
\end{align}
where stationary phase estimates as in \cite{Mockenhaupt96,CLMP} were used  in the
penultimate line. A direct computation as in \cite{CLMP} gives
\[
 \int_{\R^d}
 \left|e^{-\xi^T(I-i\overline{Q}(\eta))^{-1}\xi/2}\right|
 \,d\xi
 \sim \det(I+\overline{Q}(\eta)^2)^{1/2}.
\]
We note that, for any $A\in M_{d\times d}(\R)$,
$\det(I+A^2)=|\det(I-iA)|^2$.
Thus, combining \eqref{GaussianSplit} with the above
simplification yields
\[
 \int_{\R^d}
 \left|\widehat{e^{-|x|^2/2}\,d\sigma}(\xi,\eta)\right|
 \,d\xi
 \sim \det(I+\overline{Q}(\eta)^2)^{1/4}.
\]
In order to pass the estimate to a compactly supported smooth
density $a\,d\sigma$, we may write
$a\,d\sigma=\widetilde a e^{-|x|^2/2}\,d\sigma$,
where $\widetilde a(x)=a(x)e^{|x|^2/2}$.
For a fixed $\eta$, using Young's convolution inequality,
we obtain
\[
 \int_{\R^d}|\widehat{a\,d\sigma}(\xi,\eta)|\,d\xi
 \lesssim_a
 \int_{\R^d}
 \left|\widehat{e^{-|x|^2/2}\,d\sigma}(\xi,\eta)\right|
 \,d\xi.
\]
We estimate the $L^1$ norm of $\widehat{ad\sigma}$ over an annular neighborhood
of the normal directions with slower decay. For $R\ge2$, let $A_{\mathbf a}(R)\subset\S^{n-d-1}$ be an angular
region, to be specified below, and set
\[
 \Omega_{\mathbf a}(R)
 :=
 \left\{(\xi,\eta)\in\R^d\times\R^{n-d}:
 \frac R2\le|\eta|<R,\quad
 \frac{\eta}{|\eta|}\in A_{\mathbf a}(R)\right\}.
\]
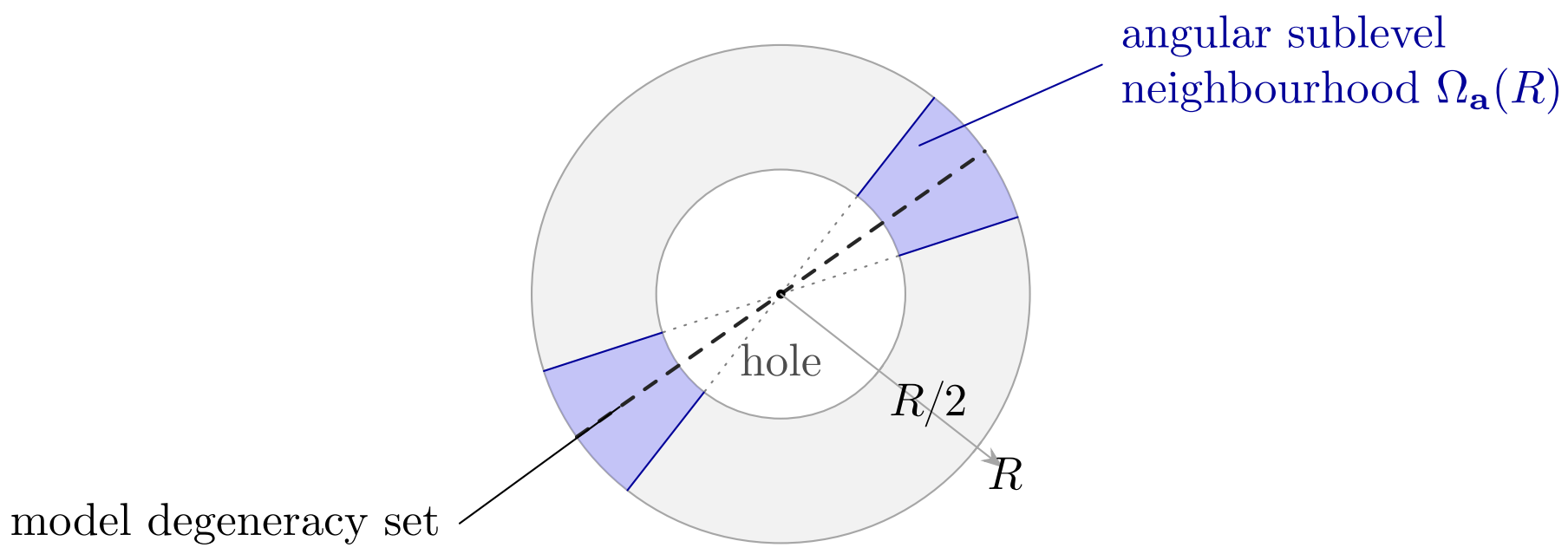
\begin{figure}[ht]
\centering
\begin{tikzpicture}[>=Stealth,line cap=round]
\def\Ro{1.9}\def\Ri{0.95}\def\th{35}\def\dl{17}
\fill[gray!10,even odd rule] (0,0) circle (\Ro) (0,0) circle (\Ri);
\draw[gray!70,thin] (0,0) circle (\Ro);
\draw[gray!70,thin] (0,0) circle (\Ri);
\foreach \b in {\th,{\th+180}}{
  \pgfmathsetmacro{\lo}{\b-\dl}\pgfmathsetmacro{\hi}{\b+\dl}
  \fill[blue!45,fill opacity=0.45]
    (\lo:\Ri) -- (\lo:\Ro) arc[start angle=\lo,end angle=\hi,radius=\Ro]
    -- (\hi:\Ri) arc[start angle=\hi,end angle=\lo,radius=\Ri] -- cycle;
  \draw[blue!60!black,thin] (\lo:\Ri) -- (\lo:\Ro) (\hi:\Ri) -- (\hi:\Ro);
  \draw[blue!60!black,gray,dotted] (0,0) -- (\lo:\Ri) (0,0) -- (\hi:\Ri);
}
\draw[black!85,dashed,thick] ({-\Ro*cos(\th)},{-\Ro*sin(\th)}) -- ({\Ro*cos(\th)},{\Ro*sin(\th)});
\fill (0,0) circle (1pt);
\node[font=\footnotesize,gray!60!black] at (0,-0.5) {hole};
\draw[->,gray!70,thin] (0,0) -- (-38:{\Ro+0.25});
\foreach \r/\l in {\Ri/{R/2},\Ro/{R}}{
  \draw[gray!70,thin] (-38:\r) ++(-38+90:0.05) -- ++(-38-90:0.1);
  \node[font=\footnotesize,anchor=north west,inner sep=2pt] at (-38:\r) {$\l$};
}
\node[font=\footnotesize,blue!60!black,align=left,anchor=west] (A) at (2.45,1.75) {angular sublevel\\ neighbourhood $\Omega_{\mathbf a}(R)$};
\draw[blue!60!black,thin] (A.west) -- (47:1.55);
\node[font=\footnotesize,align=right,anchor=east] (B) at (-2.45,-1.75) {model degeneracy set};
\draw[thin] (B.east) -- (215:1.5);
\end{tikzpicture}
\caption{Schematic slice of $\Omega_{\mathbf a}(R)$ in the normal-frequency plane. The light annulus is
$R/2\le|\eta|<R$, with the central hole $|\eta|<R/2$ excluded. The dashed line is a model degeneracy set
(a line through the origin), and the blue region is a thickened angular neighbourhood of its two opposite radial
segments inside the annulus.}
\label{fig:Omega}
\end{figure}
Combining the above bounds, we estimate the growth of the $L^1$ norm on this region as
\[
\begin{aligned}
 \|\widehat{a\,d\sigma}\|_{L^1(\Omega_{\mathbf a}(R))}
 &\lesssim_a
 \int_{\substack{R/2\le|\eta|<R\\
                 \eta/|\eta|\in A_{\mathbf a}(R)}}
 \left|\det(I-i\overline Q(\eta))\right|^{1/2}\,d\eta\\
 &\lesssim_{a}
 \int_{\substack{R/2\le|\eta|<R\\
                 \eta/|\eta|\in A_{\mathbf a}(R)}}
 \prod_{i=1}^{d}
 \left(1+|\eta|
 \left|\mu_i\left(\frac{\eta}{|\eta|}\right)\right|
 \right)^{1/2}\,d\eta,
\end{aligned}
\]
where $\mu_i$ are the eigenvalues of the normalized Hessian. We perform a polar change of variables and let $\Gamma_i(\theta)$
denote $|\mu_i(\theta)|$ arranged in descending order.
Since the maximum rank of the Hessian is $k$, we have
$\Gamma_{k+1}\equiv\Gamma_{k+2}\equiv\cdots\equiv\Gamma_d\equiv0$,
with this condition vacuous when $k=d$. This implies
\[
 \|\widehat{a\,d\sigma}\|_{L^1(\Omega_{\mathbf a}(R))}
 \lesssim_a
 \int_{R/2}^{R}\int_{A_{\mathbf a}(R)}
 r^{n-d-1}\prod_{i=1}^{k}
 (1+r\Gamma_i(\theta))^{1/2}\,d\theta\,dr.
\]

We may expand the product as
\[
 \prod_{i=1}^{k}(1+r\Gamma_i(\theta))
 =1+\gamma_1(\theta)r+\cdots+\gamma_k(\theta)r^k.
\]
Since the minimum rank of the Hessian is $k_0$, continuity of eigen-values
and compactness of $\S^{n-d-1}$ give
\[
 \Gamma_i(\theta)\sim1,
 \qquad i\in\{1,\ldots,k_0\}.
\]
Further, choose $\mathbf a=(a_{k_0+1},\ldots,a_k)\in(0,1]^{k-k_0}$
and take
\[
 A_{\mathbf a}(R)
 :=
 \left\{\theta\in\S^{n-d-1}:
 \Gamma_i(\theta)\le C_0R^{-a_i}
 \text{ for }i\in\{k_0+1,\ldots,k\}\right\},
\]
where $C_0>0$ is fixed independently of $R$.

Then, for $\theta\in A_{\mathbf a}(R)$ and $r\ge1$, the product
is bounded above by
\[
 \prod_{i=1}^{k}(1+r\Gamma_i(\theta))
 \lesssim
 r^{k_0}
 \left(
  1+\sum_{j=1}^{k-k_0}
  r^jR^{-\sum_{i=k_0+1}^{k_0+j}a_i}
 \right).
\]
Here we used the descending ordering of the $\Gamma_i$: every
product of $j$ distinct eigenvalues among
$\Gamma_{1},\ldots,\Gamma_k$ is bounded above by
$\Gamma_{1}\cdots\Gamma_{j}$.

Consequently, for $R\ge2$,
\[
\begin{aligned}
 \|\widehat{a\,d\sigma}\|_{L^1(\Omega_{\mathbf a}(R))}
 &\lesssim_a
 |A_{\mathbf a}(R)|
 \int_{R/2}^{R}
 r^{n-d-1+k_0/2}
 \left(
  1+\sum_{j=1}^{k-k_0}
  r^jR^{-\sum_{i=k_0+1}^{k_0+j}a_i}
 \right)^{1/2}\,dr\\
 &\lesssim_a
 |A_{\mathbf a}(R)|
 R^{n-d+\frac{k}{2}
       -\frac12\sum_{i=k_0+1}^{k_0+j}a_i}.
\end{aligned}
\]
In the last inequality, we used $r\sim R$ and $0\le a_i\le1$.

To estimate the corresponding $L^2$ norm, we apply Plancherel
in the tangential frequencies at each fixed $\eta$. This yields,
\[
 \int_{\R^d}
 |\widehat{a\,d\sigma}(\xi,\eta)|^2\,d\xi
 =C_a,
\]
where $C_a>0$ is independent of $\eta$ for a nonzero density.
It follows that
\[
 \|\widehat{a\,d\sigma}\|_{L^2(\Omega_{\mathbf a}(R))}^2
 \approx_{a,n}
\,
 |A_{\mathbf a}(R)|R^{n-d}.
\]

Hence, whenever $|A_{\mathbf a}(R)|>0$, 
\[
 \boxed{
 \frac{\|\widehat{a\,d\sigma}\|_{L^1(\Omega_{\mathbf a}(R))}}
      {\|\widehat{a\,d\sigma}\|_{L^2(\Omega_{\mathbf a}(R))}}
 \lesssim_a
 |A_{\mathbf a}(R)|^{1/2}
 R^{\frac12\left(
 n-d+k-\sum_{i=k_0+1}^{k}a_i
 \right)}.
 }
\]
Equivalently, with the ambient normalization,
\begin{equation}
\label{eq:QuadFR}
 R^{-n/2}
 \frac{\|\widehat{a\,d\sigma}\|_{L^1(\Omega_{\mathbf a}(R))}}
      {\|\widehat{a\,d\sigma}\|_{L^2(\Omega_{\mathbf a}(R))}}
 \lesssim_a
 |A_{\mathbf a}(R)|^{1/2}
 R^{-\frac12\left(
 d-k+\sum_{i=k_0+1}^{k}a_i
 \right)}.
\end{equation}
We now use the Fourier-ratio bound to obtain a lower bound on
$L^p$ growth and a necessary condition for $L^p$ integrability.

\textbf{Step II: Sublevel-set estimates and $L^p$ growth}:\\
For $p>2$, interpolation yields
\[
 \|\widehat{a\,d\sigma}\|_{L^2(\Omega_{\mathbf a}(R))}
 \le
 \|\widehat{a\,d\sigma}\|_{L^1(\Omega_{\mathbf a}(R))}^{\frac{p-2}{2(p-1)}}
 \|\widehat{a\,d\sigma}\|_{L^p(\Omega_{\mathbf a}(R))}^{\frac{p}{2(p-1)}}.
\]
Rearranging gives
\[
 \|\widehat{a\,d\sigma}\|_{L^p(\Omega_{\mathbf a}(R))}^p
 \ge
 \frac{\|\widehat{a\,d\sigma}\|_{L^2(\Omega_{\mathbf a}(R))}^p}
 {\left(
 \|\widehat{a\,d\sigma}\|_{L^1(\Omega_{\mathbf a}(R))}/
 \|\widehat{a\,d\sigma}\|_{L^2(\Omega_{\mathbf a}(R))}
 \right)^{p-2}}.
\]
The previously established Fourier-ratio estimate followed by
the $L^2$ estimate via Plancherel implies
\[
\|\widehat{a\,d\sigma}\|_{L^p(\Omega_{\mathbf a}(R))}^p
 \gtrsim_{a,p}
 |A_{\mathbf a}(R)|
 R^{n-d-K(\mathbf a)(p/2-1)}\]
 where $K(\textbf a)=k- \sum_{i=k_{0}+1}^{k}a_{i}$. Using sub-level set estimate,
 \[|A_{\textbf a}(R)| \gtrsim_{a} R^{-\beta_{1}(a)}\]
 in above yields,
 \[\|\widehat{a\,d\sigma}\|_{L^p(\Omega_{\mathbf a}(R))}^p
 \gtrsim_{a,p}
 R^{n-d-\beta_{1}(\mathbf a)-K(\mathbf a)(p/2-1)}.
\]
Since $\Omega_{\mathbf a}(R)\subset\{|\eta|\ge R/2\}$, we have
\[
 \|\widehat{a\,d\sigma}\|_{L^p(\{|\eta|\ge R/2\})}^p
 \gtrsim_{a,p}
 R^{n-d-\beta_{1}(\mathbf a)-K(\mathbf a)(p/2-1)}.
\]
If $\widehat{a\,d\sigma}\in L^p$, the left-hand side tends to
zero as $R\to\infty$.
Thus, provided $K(\mathbf a)>0$,
$\widehat{a\,d\sigma}\notin L^p$ whenever
\[
 2<p\le
 \frac{2(n-d+K(\mathbf a)-\beta_{1}(\mathbf a))}{K(\mathbf a)}.
\]
We now produce an upper bound for $p_{c}$.

\textbf{Step-III: $L^{p}$  upper bound via Dyadic slicing}\\
Following the Gaussian computation in Step~I and using Young's
convolution inequality, we obtain
\[
 \int_{\mathbb R^d}
 |\widehat{a\,d\sigma}(\xi,r\theta)|^p\,d\xi
 \lesssim_{a,p}
 \prod_{i=1}^{k}(1+r\Gamma_i(\theta))^{1-p/2}.
\]
For a more detailed calculation of the above, an interested reader may
refer to \cite{CLMP}.
Put $q=p/2-1>0$. Choose $C_0>0$ so that
$\Gamma_i(\theta)\le C_0$ for all $i$ and $\theta$.
For $R\ge2$ and $0\le u\le C_0$, we may use
\[
 \max\{R^{-1},u/C_0\}^{-q}
 =
 1+q\int_{R^{-1}}^1
 \mathbf 1_{\{u\le C_0t\}}t^{-q-1}\,dt
\]
to obtain
\[
 (1+Ru)^{-q}
 \lesssim_{q,C_0}
 R^{-q}
 \left(
  1+q\int_{R^{-1}}^1
  \mathbf 1_{\{u\le C_0t\}}t^{-q-1}\,dt
 \right).
\]
Changing variables $t=R^{-s}$, we obtain
\[
 (1+Ru)^{-q}
 \lesssim_{q,C_0}
 R^{-q}
 \left(
  1+q\log R\int_0^1
  \mathbf 1_{\{u\le C_0R^{-s}\}}R^{qs}\,ds
 \right).
\]
For $R/2\le r\le R$, we may therefore write
\[
 \prod_{i=1}^{k}(1+r\Gamma_i(\theta))^{-q}
 \lesssim_{k,q,C_0}
 R^{-kq}
 \prod_{i=1}^{k}
 \left(
  1+q\log R\int_0^1
  \mathbf 1_{\{\Gamma_i(\theta)\le C_0R^{-s}\}}
  R^{qs}\,ds
 \right).
\]
Since $\Gamma_i\sim1$ for $i\le k_0$, the corresponding
factors on the right are uniformly bounded. Expanding the
remaining product gives
\[
\begin{aligned}
 &\int_{\S^{n-d-1}}
 \prod_{i=1}^{k}(1+R\Gamma_i(\theta))^{-q}\,d\theta\\
 &\quad\lesssim_{p,\overline Q,C_0}
 R^{-qk}
 \sum_{J\subseteq\{k_0+1,\ldots,k\}}
 (q\log R)^{|J|}
 \int_{[0,1]^J}
 R^{q\sum_{i\in J}a_i}
 \left|
 \left\{\theta:
 \Gamma_i(\theta)\le C_0R^{-a_i},\ i\in J
 \right\}
 \right|\,d\mathbf a_J.
\end{aligned}
\]
The set inside the integral is $A_{\mathbf a}(R)$, as
defined earlier, with $a_i=0$ for $i\notin J$.
For $J=\varnothing$, the integral is understood to be
$|\S^{n-d-1}|$. The uniform sublevel-volume upper bound gives
\[
 |A_{\mathbf a}(R)|
 \lesssim R^{-\beta_2(\mathbf a)}.
\]
Set
\[
 p_0=\sup_{\mathbf a}
 \frac{2(n-d+K(\mathbf a)-\beta_2(\mathbf a))}
      {K(\mathbf a)}.
\]
Then, for $p>p_0$,
\[
\begin{aligned}
 n-d-\beta_2(\mathbf a)-qK(\mathbf a)
 &\le -\frac{K(\mathbf a)}2(p-p_0)\\
 &\le -\frac{k_0}2(p-p_0).
\end{aligned}
\]
Consequently, integration over the normal shell gives
\[
\begin{aligned}
 &\int_{R/2\le|\eta|<R}\int_{\mathbb R^d}
 |\widehat{a\,d\sigma}(\xi,\eta)|^p\,d\xi\,d\eta\\
 &\quad\lesssim_{a,p}
 R^{n-d-qk}
 \sum_{J\subseteq\{k_0+1,\ldots,k\}}
 (q\log R)^{|J|}
 \int_{[0,1]^J}
 R^{q\sum_{i\in J}a_i-\beta_2(\mathbf a)}
 \,d\mathbf a_J\\
 &\quad\lesssim_{a,p}
 (1+\log R)^{k-k_0}
 R^{-k_0(p-p_0)/2}.
\end{aligned}
\]
Summing over $R=2^j$ yields
\[
 \int_{|\eta|\ge1}\int_{\mathbb R^d}
 |\widehat{a\,d\sigma}(\xi,\eta)|^p\,d\xi\,d\eta
 \lesssim_{a,p}
 \sum_{j\ge1}
 (1+j)^{k-k_0}2^{-jk_0(p-p_0)/2}.
\]
Since $k_0>0$, the sum on the right converges when $p>p_0$,
and the contribution from $|\eta|<1$ is finite by the
Gaussian upper bound. Hence
$\widehat{a\,d\sigma}\in L^p(\mathbb R^{n})$
for every finite $p>p_0$. Membership in $L^\infty$
follows from the finite total variation of $a\,d\sigma$.
This implies
\[
 \boxed{
 p_c\le
 \sup_{\mathbf a}
 \frac{2(n-d+K(\mathbf a)-\beta_2(\mathbf a))}
      {K(\mathbf a)}.
 }
\]
We now extend the results to all compactly supported distributions.

\textbf{Step IV: Extension to distributions}\\
Suppose, $0\ne T$ is a
compactly supported Tempered distribution on the quadratic graph and
$\widehat T\in L^p$ for some
\[
 2\le p\le \sup_{\mathbf a}
 \frac{2(n-d+K(\mathbf a)-\beta_1(\mathbf a))}{K(\mathbf a)}.
\]
In what follows, we closely follow the proof of Lemma~0.9,
adapted to the region $\Omega_{\mathbf a}(R)$.
Define the smooth flattening map
$S:U\times\R^{n-d}\rightarrow U\times\R^{n-d}$ by
\[
 S(t',t'')=(t',t''-\Phi(t')).
\]
Let $\alpha_0$ index a transverse derivative of maximal order
with nonzero coefficient in the expansion of the flattened
distribution $S_*T$, as in Lemma~0.9. Choose $\rho\in C_c^\infty$ equal to $1$ near
$\operatorname{supp}T$, and set
\[
 T_0=\rho(t',t'')\bigl(t''-\Phi(t')\bigr)^{\alpha_0}T.
\]
The maximality of $|\alpha_0|=m$ gives
\[
 S_*T_0=v\otimes\delta_0,
 \qquad
 v=(-1)^m\alpha_0!\,T_{\alpha_0}\ne0.
\]
Moreover, $\widehat T_0\in L^p$ by Young's inequality.
We may choose $a\in C_c^\infty(U)$ with
$\langle v,a\rangle\ne0$.
For an admissible $\mathbf a$, set
\[
 D_R=
 \{\eta:R/2<|\eta|<R,\ \eta/|\eta|\in A_{\mathbf a}(R)\},
 \qquad
 M_R=|D_R| \sim R^{n-d}|A_{\mathbf a}(R)|. 
\]
Consider $h_{R} \in C_{c}^{\infty}(\R^{n-d})$, such that $\operatorname{Supp}   h_{R} \subseteq D_{R}$, $0\leq h_{R} \leq 1$. Since $D_R$ is semialgebraic, its boundary has measure zero\cite{BCR98}. Thus, we may choose $h_R$ so that
it additionally satisfies $\int h_R\sim M_R$. Set $b_{R}= \hat{h}_{R}$ and define
\[
 g_R(t',t'')=a(t')b_R(t''-\Phi(t')).
\]
By the definition of $T_0$,
\[
 |\langle T_0,g_R\rangle|
 =b_R(0)|\langle v,a\rangle|
 \sim_{T_0,a}M_R,
 \qquad
 \widehat g_R(\xi,\eta)
 =h_R(\eta)\widehat{a\,d\sigma}(\xi,\eta).
\]
Using the previously established $L^1$ estimate, we obtain
\[
\begin{aligned}
 \|\widehat g_R\|_{L^1(\R^d\times D_R)}
 &=\int_{D_R}h_R(\eta)
   \int_{\R^d}|\widehat{a\,d\sigma}(\xi,\eta)|\,d\xi\,d\eta\\
 &\lesssim_a
 R^{K(\mathbf a)/2}\int_{D_R}h_R(\eta)\,d\eta
 \lesssim_a M_RR^{K(\mathbf a)/2}.
\end{aligned}
\]
Applying Plancherel's theorem yields the $L^2$ estimate
\[
 \|\widehat g_R\|_{L^2(\R^d\times D_R)}
 \lesssim_a M_R^{1/2}.
\]
Combining the above two estimates by interpolation, with
$p'=p/(p-1)$, yields
\[
 \|\widehat g_R\|_{p'}
 \lesssim_{a,p}
 M_R^{1-1/p}R^{K(\mathbf a)(1/2-1/p)}.
\]
Fourier duality and H\"older's inequality therefore imply
\[
\begin{aligned}
 \|\widehat T_0\|_{L^p(\{|\eta|\ge R/2\})}
 &\gtrsim_{T_0,a,p}
 M_R^{1/p}R^{-K(\mathbf a)(1/2-1/p)}\\
 &\gtrsim_{T_0,a,p}
 R^{\frac{n-d-\beta_1(\mathbf a)}p
       -K(\mathbf a)(\frac12-\frac1p)}.
\end{aligned}
\]
For any admissible $\mathbf a$ with $K(\mathbf a)>0$, if
\[
 p\le
 \frac{2(n-d+K(\mathbf a)-\beta_1(\mathbf a))}{K(\mathbf a)},
\]
the lower bound on the right-hand side forces $\widehat T_0$
to have a nonvanishing $L^p$ tail. This contradicts
$\widehat T_0\in L^p(\R^{n})$.
Taking the supremum over admissible $\mathbf a$ and combining
this with the conclusion of the previous step yields
\[
 \boxed{
 \sup_{\mathbf a}
 \frac{2(n-d+K(\mathbf a)-\beta_1(\mathbf a))}
      {K(\mathbf a)}
 \leq p_c\le
 \sup_{\mathbf a}
 \frac{2(n-d+K(\mathbf a)-\beta_2(\mathbf a))}
      {K(\mathbf a)}.
 }
\]
\end{proof}

\subsection{ Proof of Theorem\ref{thm:PthinConstantNullity}}
\begin{proof}

\textbf{Step I: Fourier ratio and lower bound}

We first prove an $L^{p}$ tail bound for the smooth density and then pass the estimate to a general distribution. Suppose we locally write $M$ as the graph
\[
M_0=\{(x,\Psi(x))~|~x\in U\subseteq \R^{d}\},
\]
where $U$ is open and $\Psi:U\to\R^{n-d}$ is a smooth function parametrizing the manifold of constant relative nullity, with $\nabla\Psi(0)=0$. Let $\sigma_M$ denote the surface measure on $M$ near the origin, and let $0 \neq f\in C_c^\infty(M)$. We then use M\"uller's theorem\cite{Muller82} to bound the tangential $L^{1}$ norm:
\[
\int_{\R^d}\left|\widehat{f\,d\sigma_M}(\xi,\eta)\right|\,d\xi
\lesssim_{f,d,n} (1+|\eta|)^{\frac{\rho}{2}},
\qquad \eta\in\R^{n-d}.
\] Put
\[
\Omega_R=\mathbb R^d\times\{\eta:R<|\eta|<2R\}.
\]
Then, using the tangential $L^{1}$ bound, we get
\[
\|\widehat{f\,d\sigma}\|_{L^1(\Omega_R)}
\lesssim_f R^{c+\rho/2}.
\]
Writing $d\sigma=J_\Psi(x)\,dx$ in graph coordinates, we first use tangential Plancherel to get
\[
\int_{\R^d}
\left|\widehat{f\,d\sigma}(\xi,\eta)\right|^2\,d\xi
=
\int_U |f(x,\Psi(x))|^2 J_\Psi(x)^2\,dx.
\]
We then integrate over $R<|\eta|<2R$ to obtain, 
\[
\|\widehat{f\,d\sigma}\|_{L^2(\Omega_R)}^2
=
\int_{R<|\eta|<2R}
\int_U |f(x,\Psi(x))|^2 J_\Psi(x)^2\,dx\,d\eta
\sim_f R^c.
\]
Thus the annular Fourier ratio has decay
\begin{equation}
R^{-n/2}
\frac{\|\widehat{f\,d\sigma}\|_{L^1(\Omega_R)}}
{\|\widehat{f\,d\sigma}\|_{L^2(\Omega_R)}}
\lesssim_f R^{-\nu/2}.
\end{equation}
Interpolation between $L^1$ and $L^p$ as in\cite{DI26} yields, for $p>2$,
\begin{equation}
\|\widehat{f\,d\sigma}\|_{L^p(\Omega_R)}^p
\ge
\frac{\|\widehat{f\,d\sigma}\|_{L^2(\Omega_R)}^{2p-2}}
{\|\widehat{f\,d\sigma}\|_{L^1(\Omega_R)}^{p-2}}
\gtrsim_{f,p} R^{((n-d)-\rho(p/2-1)}.
\end{equation}
As we have obtained $L^{p}$ tail bounds for smooth densities supported on $M$, we follow the proof of Lemma~\ref{lem:annular-testing} to obtain, for every nonzero compactly supported distribution $T$ on $M$ and $2\le p<\infty$,
\begin{equation}\label{eq:NullityTail}
    \|\hat{T}\|_{L^{p}(\Omega(R))}^{p}
\gtrsim_{T,p} R^{((n-d)-\rho(\frac{p}{2}-1)}.
\end{equation}
If $\hat{T}\in L^{p}$, its $L^{p}$ mass on $\Omega(R)$ must tend to zero as $R\to\infty$. Thus the exponent on the right must be non-positive, which requires
\[
p\geq2+\frac{2(n-d)}{\rho}.
\]
We therefore conclude
\[
\boxed{P_c\ge 2+\frac{2(n-d)}{\rho}}.
\]

\medskip
\textbf{Step II: Pointwise decay and upper bound}

We use the pointwise decay bound obtained in Theorem~\ref{Prop:Uniformrankdecay} to establish $L^p$ membership of the Fourier transform of a smooth density for $p$ above the claimed upper bound. As in Step I, let $\sigma$ denote the surface measure on $M$, and choose a nonzero $f\in C_c^\infty(M).$ As every nonzero normal combination of the graph Hessians has rank at least \(r>0\) the uniform rank estimate from \ref{Prop:Uniformrankdecay} yield
\begin{equation}
\|\widehat{fd\sigma}(\zeta,\xi,\eta)\|_\infty\lesssim_{d,n,f} (1+|\zeta|+|\xi|+|\eta|)^{\frac{-r}{2}}.
\end{equation}
Denote by $J_\Psi$ the surface Jacobian of the local graph
parametrization $x\mapsto (x,\Psi(x))$ of $M$. Then tangential Plancherel in
$(\zeta,\xi)$ gives, for each fixed $\eta$,
\begin{equation}\label{eq:plancherel}
\begin{aligned}
\int_{\xi}
\|\widehat{f\,d\sigma}(\cdot,\xi,\eta)\|_{L^2_\zeta}^2\,d\xi
&=(2\pi)^{n}
\int \bigl|f(x,\phi(x))J_\phi(x)
           e^{ i\eta\cdot\phi(x)}\bigr|^2\,dx \\
&=(2\pi)^{n}
\|f(\cdot,\phi(\cdot))J_\phi(\cdot)\|_{L^2(dx)}^2
\lesssim_{\phi} \|f\|_{L^2(M,d\sigma)}^2,
\end{aligned}
\end{equation}
Therefore, for $p>2$, $L^{2}-L^{\infty}$ interpolation in $\zeta$ followed by integration in $(\xi,\eta)$ gives
\begin{equation}
\begin{split}
\|\widehat{f\,d\sigma}\|_{L^p(\R^n)}^p
&\lesssim_{f,p}
\int_{\R^c}\int_{\xi}
(1+|\xi|+|\eta|)^{-\frac{r(p-2)}{2}}
\|\widehat{f\,d\sigma}(\cdot,\xi,\eta)\|_{L^2_\zeta}^{2}
\,d\xi\,d\eta\\
&\le
\int_{\R^c}
(1+|\eta|)^{-\frac{r(p-2)}{2}}
\int_{\xi}
\|\widehat{f\,d\sigma}(\cdot,\xi,\eta)\|_{L^2_\zeta}^{2}
\,d\xi\,d\eta\\
&\lesssim_{f,p}
\int_{\R^c}
(1+|\eta|)^{-\frac{r(p-2)}{2}}\,d\eta.
\end{split}
\end{equation}
The last integral is finite whenever
\[
\frac{r(p-2)}{2}>n-d.
\]
Thus
\[
\widehat{f\,d\sigma}\in L^p(\R^n)
\qquad\text{for every}\qquad
p>2+\frac{2(n-d)}{r}.
\]
Combining above with lower bound at the end of step-I we conclude
\[
\boxed{2+\frac{2(n-d)}{\rho} \leq P_c\le 2+\frac{2(n-d)}{r}}.
\]
\end{proof}
\begin{remark}
\begin{remark}
In the hypersurface case, the upper and lower bounds coincide,
recovering Guo's~\cite{Guo93} sharp exponent
\[
p_c=2+\frac{2}{\rho}.
\]
This gives an alternative proof of Guo's result for hypersurfaces
with constant relative nullity, using M\"uller's theorem,
interpolation, and the annular testing argument developed
in Lemma~\ref{lem:annular-testing}.
\end{remark}
\end{remark}
\subsection{Proof Of Theorem\ref{QuadraticImprovement}}
\begin{proof}
Taking $g\equiv1$, we have $g\in L^p(M_Q,d\sigma)$ for every
$1\leq p\leq\infty$, since $\sigma$ is finite. Consider a dyadic shell
$R<|\eta|<2R$. Fix an admissible $h$ and write $\alpha=(t-2)/2$,
where $t>2$ and $\alpha\geq h$. Using the tangential $L^1$ estimate
from Step~I and the interpolation argument from Step~II of the proof
of Theorem~\ref{Thm:QuadraticPthin}, followed by the change of variables
$\eta=Rw$, gives
\[
\int_{R<|\eta|<2R}\int_{\mathbb R^d}
|E_\sigma g(\xi,\eta)|^t\,d\xi\,d\eta
\gtrsim
R^{m-d\alpha}\int_{1<|w|<2}D_{1/R}(w)^{-\alpha}\,dw.
\]
By Lemma~\ref{lem:determinant-integral}, the integral on the right is
$\gtrsim R^{h\alpha-h(h+1)/2}$ for $h\geq1$. For $h=0$, it is
$\gtrsim1$ since the singular values are bounded above on $1<|w|<2$.

We now sum over the dyadic shells $R=2^j$ to obtain
\[
\|E_\sigma g\|_{L^t(\mathbb R^{d+m})}^t
\gtrsim
\sum_{j\geq1}2^{j\left(m-(d-h)\alpha-h(h+1)/2\right)}.
\]
Taking $\alpha=(m-h(h+1)/2)/(d-h)$, which exceeds $h$ by
admissibility, makes the exponent zero. Thus the sum diverges at
\[
t=2+\frac{2m-h(h+1)}{d-h}.
\]
Since $E_\sigma g$ is bounded, this also excludes every smaller
exponent. Taking the maximum over admissible $h$ gives
\[
E_\sigma g\in L^q(\mathbb R^{d+m})
\implies q>P_{\mathrm{quad}}(d,m).
\]

Finally, linear independence gives $m\leq d(d+1)/2$, so $D=d+2m$.
The oblique inequality from $\tau_0(d,d+m)$ therefore gives
\[
\frac1p+\frac{d+2m}{dq}\leq1.
\]
Together, these inequalities yield
$\mathcal T_\sigma(M_Q)\subseteq\tau_Q(d,m)$.
\end{proof}
\subsection{Proof Of Theorem\ref{thm:PDEConstantnullity}}
\begin{proof}
    Suppose $u \in \mathcal{S}'(\mathbb{R}^n) \cap L^p_{\mathrm{loc}}(\mathbb{R}^n)$ 
satisfies the hypotheses of the theorem, and suppose moreover that $Pu = 0$. 
Taking the Fourier transform gives
\[
    P(\xi)\,\widehat{u}(\xi) = 0 \qquad \text{in } \mathcal{S}'(\mathbb{R}^n),
\]
from which it follows that $\operatorname{supp}\widehat{u}$ is contained in the 
zero set of $P$. We endow this zero set with a locally finite 
stratification\cite{BCR98}
\[
    \{\xi \in \mathbb{R}^n : P(\xi) = 0\} = \bigcup_{j} V_j,
\]
where each $V_j$ is a smooth manifold, the strata are pairwise disjoint, and 
 $V_i \subset \overline{V_j}$ whenever $i < j$, 
$\dim V_i < \dim V_j$, and $V_i \cap \overline{V_j} \neq \emptyset$. Let $V_{\alpha}$ be the stratum of maximum dimension in $\operatorname{supp}\widehat{u}$. 
The set
\[
    V =V_{\alpha} \setminus \bigcup_{j < \alpha} \overline{V_j}
\]
is open in $V_{\alpha}$, and hence $V = U \cap V_{\alpha}$ for some open set 
$U \subset \mathbb{R}^n$ which may be chosen disjoint from all $V_j$ with 
$j < \alpha$. Choose $\rho \in \mathcal{S}(\mathbb{R}^n)$ with
$\widehat{\rho} \in C_c^{\infty}(\mathbb{R}^n)$ and
$\operatorname{supp}\widehat{\rho} \subseteq U$, and consider $\hat{u}\cdot\hat{\rho}$, then
$$\operatorname{Supp}(\hat{u}\cdot\hat{\rho}) \subseteq \operatorname{Supp}(\hat{u})\cap \operatorname{Supp}(\hat{\rho}) \subseteq V_{j}. $$
Set $T_{\alpha}= \hat{u}\cdot\hat{\rho}$, if $T_{\alpha}$ is not zero, then by lemma \ref{lem:annular-testing},
$$\|\hat{T}_{\alpha}\|_{L^p(\Gamma_R)}
 \gtrsim_{T,p,\varepsilon}
 R^{\,m_{j}+\frac{n-2\kappa_\sigma}{p}
             -\frac{d-2\kappa_\sigma}{2}-\varepsilon},
 \qquad R\ge R_0(T,p,\varepsilon).$$
By local decay condition \eqref{eq:Stratum_Lp-decay}, we may find a sequence $\{R_{n}\}_{n} $ such that,
$$\|\hat{T}_{\alpha}\|_{L^p(\Gamma_{R_{n}})} \lesssim R_{n}^{\frac{s_{j}}{p}}.$$
Comparing the exponent $R$ from above,
$$\,m_{j}+\frac{n-2\kappa_\sigma}{p}
             -\frac{d-2\kappa_\sigma}{2}-\varepsilon \leq \frac{s_{j}}{p} .$$
The bound on $p$ is choosen precisely such that,
$$\,m_{j}+\frac{n-2\kappa_\sigma}{p}
             -\frac{d-2\kappa_\sigma}{2}-\varepsilon > \frac{s_{j}}{p}.$$
This leads to contradiction, thus $T_{\alpha}=0.$ This implies $u*\rho \equiv 0$. Since the choice 
of $\rho$ was arbitrary, we conclude $\operatorname{supp}\widehat{u} \cap 
V_\alpha = \emptyset$. Proceeding inductively through the strata of strictly 
smaller dimension ($d\geq 1)$ and applying the same argument at each stage, we obtain 
$\operatorname{supp}\widehat{u} \cap (\bigcup_{j \geq 1} V_{j}) = \emptyset$. We briefly treat the case, If $V_j=\{\xi_j\}$ is a zero-dimensional stratum, every nonzero
distribution $T$ supported on $V_j$ has the form
\[
 T=\sum_{|\alpha|\le m}c_\alpha\partial^\alpha\delta_{\xi_j},
 \qquad \max_{c_\alpha\ne0}|\alpha|=m.
\]
Consequently,
\[
 \widehat T(x)=e^{ix\cdot\xi_j}Q(x),
\]
where $Q$ is a polynomial of degree $m$. Writing $Q_m$ for its
nonzero homogeneous part of degree $m$, a change of variables gives
\[
 R^{-m-n/p}\|\widehat T\|_{L^p(\Gamma_R)}
 =\|R^{-m}Q(R\,\cdot)\|_{L^p(\Gamma_1)}
 \longrightarrow \|Q_m\|_{L^p(\Gamma_1)}>0.
\]
Thus the local decay condition with $s_j\le n$ is impossible
for every finite $p$, and the localized distribution must vanish. and therefore $\operatorname{Supp}\hat{u} \cap V_{0}= 0$, thus we conclude $u \equiv0$.
\end{proof}
\begin{remark}
To deduce Theorem~\ref{thm:PDEUniqueness}, we take $m_j=0$
and use the annular decay of $u$ to estimate each localization
$u_j=u*\check\chi_j$. The proof of Lemma~\ref{lem:annular-testing}
gives the same lower bound on the smaller annulus
\[
\Gamma'_R=\{x:5R/4<|x|<7R/4\}.
\]
Choose a smooth cutoff supported in $\Gamma_R$ and equal to one
on a neighborhood of $\Gamma'_R$ of width comparable to $R$.
By Young's inequality for the part of $u$ inside the cutoff,
and by the rapid decay of $\check\chi_j$ and the temperedness
of $u$ for the remaining part, we obtain
\[
\|u_j\|_{L^p(\Gamma'_R)}
\lesssim_{\chi_j}\|u\|_{L^p(\Gamma_R)}+O_N(R^{-N})
\]
for every $N>0$. Since the testing lower bound holds for all
sufficiently large $R$, we may apply it along the same sequence
of radii supplied by \eqref{eq:Lp-decay}.
\end{remark}
\section{Proofs of Lemmas and Corollary}
In this section, we give the proofs of several lemmas,
propositions, and corollaries. The lemmas proved here
play a key role in the proofs of the  theorems and
may be useful in other settings.

\textbf{Proof of Corollary\ref{cor:cm-failure-codimension-two}}
\begin{proof}
The proof follows the same argument as that of
Theorem~\ref{thm:smooth-high-codimension}, choosing the
localization centers in the support of the reduced
tangential distribution, where CM fails.
In Step IV, we apply
Lemma~\ref{lem:cm-failure-quadratic-factor}
to obtain the required lower bound on
\[
\int_{R<|\eta|<2R}A_j(\eta)^{-(p-2)}\,d\eta.
\]
Proceeding with the same computation gives the stated
exponent, including the endpoint.
\end{proof}
We present the lemma needed for the proof of above Corollary.
\begin{lemma}[CM failure and the quadratic factor]
\label{lem:cm-failure-quadratic-factor}
Let $U\subset\mathbb R^d$ be open and let
$\Phi=(\phi_1,\phi_2)\in C^\infty(U;\mathbb R^2)$.
Suppose $K\subsetneq  U$ and CM condition fails at every point of $K$. For $x_j\in K$, write
\[
H_j(\omega)
=
\omega_1D^2\phi_1(x_j)+\omega_2D^2\phi_2(x_j), 
\]
and let $\Gamma_{j,1}(\omega)\geq \ldots\geq \Gamma_{j,d}(\omega)\geq0$
be absolute values of it's singular values. For $R>0$ and $0<\tau\leq1$, define
\[
A_j(\eta)
=
\prod_{i=1}^d
\bigl(1+\tau^{-1}\Gamma_{j,i}(\eta/R)\bigr)^{1/2}.
\]
Then for every $p>2$,
\[
\int_{R<|\eta|<2R}A_j(\eta)^{-(p-2)}\,d\eta
\gtrsim_{\Phi,K,d,p}
R^2\tau^{\frac1m+\frac{(d-1)(p-2)}2}
\]
where, $m=\lfloor\frac{d}{2}\rfloor+1$.
\end{lemma}

\begin{proof}
Fix $x_j\in K$. We expand $A_j(\eta)^2$ as
\[
A_j(\eta)^2
=
\sum_{k=0}^{d-1}\tau^{-k}
\sum_{\substack{I\subset\{1,\ldots,d\}\\ |I|=k}}
\prod_{i\in I}\Gamma_{j,i}(\eta/R)
+
\tau^{-d}\prod_{i=1}^d\Gamma_{j,i}(\eta/R).
\]
Since the singular values are uniformly bounded for
$x_j\in K$ and $R<|\eta|<2R$, all terms of degree at most
$d-1$ are bounded by a constant multiple of
$\tau^{-(d-1)}$. Thus,
\[
A_j(\eta)^2
\lesssim_{\Phi,K,d}
\tau^{-(d-1)}
+
\tau^{-d}\bigl|\det H_j(\eta/R)\bigr|.
\]
Now make the polar change of variables
$\eta=Rr(\cos\theta,\sin\theta)$, where $1<r<2$ and $\theta\in \S^1$.
By homogeneity,
\[
\Gamma_{j,i}(\eta/R)=r\Gamma_{j,i}(\theta).
\] 
and hence
\[
\bigl|\det H_j(\eta/R)\bigr|
=
r^d\prod_{i=1}^d\Gamma_{j,i}(\theta)
=
r^d\bigl|\det H_j(\theta)\bigr|.
\]
Set $P_j(\theta)=\det H_j(\theta)$. Since $1<r<2$,
we obtain
\begin{equation}{\label{eq:CMHessianBound}}
A_j(\eta)^2
\lesssim_{\Phi,K,d}
\tau^{-(d-1)}+\tau^{-d}|P_j(\theta)|.
\end{equation}
We now use CM failure to estimate $P_j$ near an
appropriate angular zero. By \cite{Christ82}, in codimension-two manifolds, failure
of the CM condition is equivalent to having a normal
direction where the determinant of the Hessian pencil
has a zero of multiplicity greater than $d/2$. Let $\theta_0$ denote one such direction, with multiplicity 
$b\geq\lfloor d/2\rfloor+1$.
This implies, in a local angular coordinate centered at $\theta_0$,
the derivatives of $P_j$ through order $m-1$ vanish.
Taylor's theorem and the uniform bounds on its angular
derivatives therefore give
\[
|P_j(\theta)|
\lesssim_{\Phi,K,d}
\operatorname{dist}_{\S^1}(\theta,\theta_0)^m
\]
for $\theta$ near $\theta_0$.
Setting
\[
E_{j,\tau}
=
\{\theta\in \S^1:
\operatorname{dist}_{\S^1}(\theta,\theta_j)<\tau^{1/m}\},
\]
we obtain $|P_j(\theta)|\lesssim\tau$ on $E_{j,\tau}$. Thus for $\theta_{0} \in E_{j,\tau}$ from \eqref{eq:CMHessianBound} we have,
\[A_{j}(Rr\theta)^{2} \lesssim_{\Phi,K,d} r^{d}\tau^{-(d-1)}.\]
Moreover, We may use arc-length formula, followed by multiplication by the radial polar-coordinate factor to get,
\[|E_{\tau,\theta_0}| \approx \tau^{\frac{1}{m}}.\]
As $p\geq2$, the above two observations yeild,
\[
\begin{aligned}
\int_{R<|\eta|<2R}A_j(\eta)^{-(p-2)}\,d\eta
&\geq
\int_R^{2R}\int_{E_{\tau,\theta_0}}
A_j(r\theta)^{-(p-2)}r\,d\theta\,dr\\
&\gtrsim
R^2\tau^{(d-1)(p-2)/2}|E_{\tau,\theta_0}|\\
&\gtrsim
R^2\tau^{\frac{(d-1)(p-2)}2+\frac1m}.
\end{aligned}
\]
\end{proof}
Next, we state the lemma required for the proof of Theorem\ref{thm:smooth-high-codimension}
\begin{lemma}\label{lem:determinant-integral}
Let
\[
 H:\mathbb R^c\longrightarrow\operatorname{Sym}_d(\mathbb R)
\]
be a bounded linear operator. For $w\in\mathbb R^c$, let
$\Gamma_1(w)\ge\cdots\ge\Gamma_d(w)\ge0$ denote the absolute values of singular
values of $H(w)$. For $0<\tau\le1$, define
\[
 D_\tau(w)=\prod_{i=1}^{d}\max\{\tau,\Gamma_i(w)\}.
\]
Suppose $1\le k\le d$ and
\[
 c>kd-\frac{k(k-1)}2.
\]
Then, for every $s\ge k$,
\[
 \int_{1<|w|<2}D_\tau(w)^{-s}\,dw
 \gtrsim_{c,k,H,s}
 \tau^{-ks+k(k+1)/2}.
\]
\end{lemma}

\begin{proof}
For fixed $w$, $D_\tau(w)^{-1}$ is related to the volume
of $\{v:|v|\le1,\ |H(w)v|\le\tau\}$. To estimate this,
we consider the set $\{(w,v):|H(w)v|\le\tau\}$, counting
its size in two ways: fixing $w$ and varying $v$ relates
its size to $D_\tau(w)$, while fixing $v$ and varying
$w$ over a suitable subset gives a lower bound.
Considering $k$ vectors simultaneously brings in the
dimensional assumption on $c$.

Fix a $k$-dimensional subspace $E\subset\mathbb R^d$.
Consider the space of vectors in $\mathbb R^c$ for which
$H(w)$ annihilates every vector in $E$:
\[
 K_E=\{w\in\mathbb R^c:H(w)v=0\text{ for every }v\in E\}.
\]
This is the kernel of the linear map
$w\mapsto H(w)|_E$ from $\mathbb R^c$ to
$\operatorname{Hom}(E,\mathbb R^d)$. By symmetry of $H(w)$
and rank--nullity,
\[
 \dim K_E
 \ge c-\bigl(d+(d-1)+\cdots+(d-k+1)\bigr)
 =c-kd+\frac{k(k-1)}2>0,
\]
where the last inequality follows from the assumption on $c$.
Let
\[
 \mathcal T_E(\tau)=
 \left\{\omega=u+z:
 u\in K_E,\ \frac54<|u|<\frac32,\quad
 z\in K_E^\perp,\ |z|<\varepsilon\tau
 \right\},
\]
where $\varepsilon>0$ will be specified later.
This is an annular thickening of the kernel $K_E$.

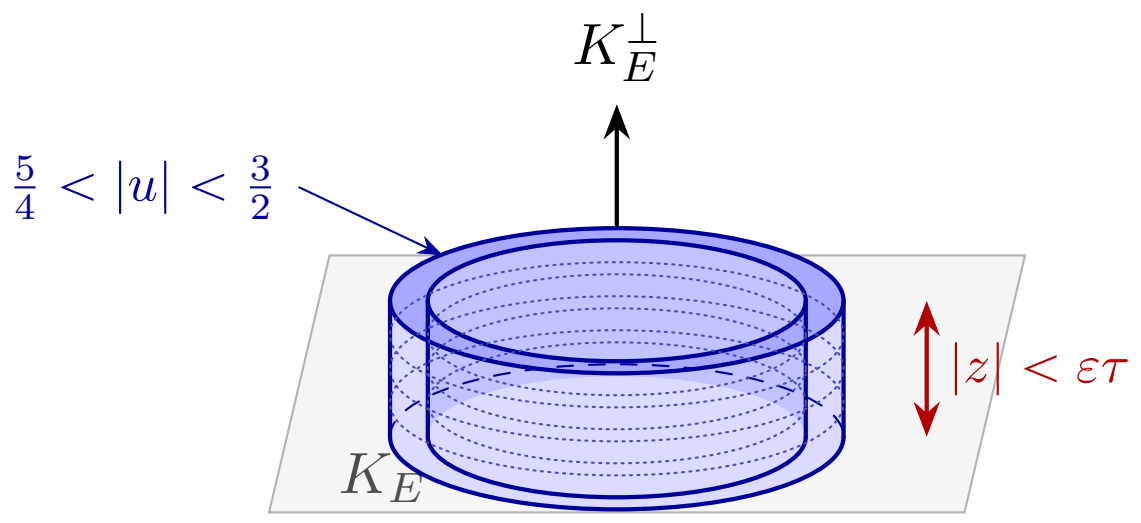
\begin{figure}[ht]
\centering
\begin{tikzpicture}[>=Stealth, line cap=round]
\def\Ro{1.5}\def\Ri{1.25}\def\sq{0.32}\def\h{0.45}
\def\bl{blue!60!black}

\fill[gray!8] (-2.3,-0.95) -- (2.3,-0.95) -- (2.7,0.75) -- (-1.9,0.75) -- cycle;
\draw[gray!60] (-2.3,-0.95) -- (2.3,-0.95) -- (2.7,0.75) -- (-1.9,0.75) -- cycle;
\node[gray!60!black,font=\small] at (-1.55,-0.72) {$K_E$};

\fill[blue!14] (-\Ro,-\h) arc (180:360:{\Ro} and {\Ro*\sq}) -- (\Ro,\h)
      arc (0:180:{\Ro} and {\Ro*\sq}) -- cycle;
\fill[white] (0,\h) ellipse ({\Ri} and {\Ri*\sq});
\fill[blue!24] (-\Ri,\h) arc (180:0:{\Ri} and {\Ri*\sq}) -- (\Ri,-\h)
      arc (0:180:{\Ri} and {\Ri*\sq}) -- cycle;
\fill[blue!34,even odd rule]
  (0,\h) ellipse ({\Ro} and {\Ro*\sq}) (0,\h) ellipse ({\Ri} and {\Ri*\sq});

\draw[\bl,thick] (-\Ro,-\h) arc (180:360:{\Ro} and {\Ro*\sq});
\draw[\bl,dashed] (-\Ro,-\h) arc (180:0:{\Ro} and {\Ro*\sq});
\draw[\bl,thick] (-\Ro,-\h) -- (-\Ro,\h) (\Ro,-\h) -- (\Ro,\h);
\draw[\bl,thick] (-\Ri,\h) -- (-\Ri,-\h) (\Ri,\h) -- (\Ri,-\h);
\draw[\bl,thick] (0,\h) ellipse ({\Ro} and {\Ro*\sq});
\draw[\bl,thick] (0,\h) ellipse ({\Ri} and {\Ri*\sq});
\draw[\bl,thick] (\Ri,-\h) arc (0:-180:{\Ri} and {\Ri*\sq});

\foreach \z in {-0.225,0,0.225}{
  \draw[\bl!70,densely dotted] (0,\z) ellipse ({\Ro} and {\Ro*\sq});
  \draw[\bl!70,densely dotted] (0,\z) ellipse ({\Ri} and {\Ri*\sq});
}

\draw[->,thick] (0,\h+0.5) -- (0,1.75) node[above,font=\small] {$K_E^{\perp}$};
\draw[red!70!black,<->,thick] (2.05,-\h) -- (2.05,\h);
\node[red!70!black,font=\footnotesize,right] at (2.05,0) {$|z|<\varepsilon\tau$};

\draw[->,\bl,thin] (-2.1,1.2) node[left,font=\footnotesize,\bl] {$\tfrac54<|u|<\tfrac32$}
   -- (-1.15,\h+0.3);
\end{tikzpicture}
\caption{The set $\mathcal T_E(\tau)$: the orthogonal thickening
of the annulus $\{u\in K_E:\tfrac54<|u|<\tfrac32\}$ by
$z\in K_E^\perp$ with $|z|<\varepsilon\tau$.}
\end{figure}

Then
\[
 |\mathcal T_E(\tau)|
 \gtrsim\tau^{\operatorname{codim}K_E}
 \ge\tau^{N_k},
\]
where $N_k=\sum_{i=0}^{k-1}(d-i)$.
Since $H$ is a bounded linear operator, we may choose
$\varepsilon>0$ sufficiently small that
\[
 \mathcal T_E(\tau)
 \subseteq
 \left\{
 \omega\in\mathbb R^c:
 1<|\omega|<2,\quad
 |H(\omega)v|\le\tau
 \text{ for every }v\in E\text{ with }|v|\le2
 \right\}.
\]
We now average over tangent directions and count these
incidences in the opposite order. Let $e_1,\ldots,e_d$ be
the standard basis vectors of $\mathbb R^d$, and choose
small fixed balls
\[
 B_i=B(e_i,\varepsilon_0),\qquad 1\le i\le k,
\]
such that $B_i\subset B(0,2)$ and every choice $v_i\in B_i$
satisfies
\[
 \left|\sum_{i=1}^{k}a_i v_i\right|
 \ge\frac12|a|,
 \qquad a\in\mathbb R^k.
\]
For each such choice, apply the above estimates with
$E=\operatorname{span}(v_1,\dots,v_k)$ to obtain
\[
 \left|
 \left\{
 w\in\mathbb R^c:
 1<|w|<2,\quad
 |H(w)v_i|\le\tau\text{ for every }i
 \right\}
 \right|
 \gtrsim\tau^{N_k}.
\]
Define
\[
 F(w)=
 \prod_{i=1}^{k}
 \left|\{v\in B_i:|H(w)v|\le\tau\}\right|.
\]
Integrating over $B_1\times\cdots\times B_k$ and applying
Fubini's theorem, we obtain
\begin{equation}\label{eq:incidence-average-lower}
 \int_{1<|w|<2}F(w)\,dw\gtrsim\tau^{N_k}.
\end{equation}

Fix $w$, and let $f_1,\ldots,f_d$ be an orthonormal basis
of eigenvectors of $H(w)$, ordered to match the singular
values $\Gamma_1(w),\ldots,\Gamma_d(w)$.
Writing $v=\sum_{j=1}^d z_j f_j$, we have
\[
 |H(w)v|^2
 =\sum_{j=1}^d\Gamma_j(w)^2|z_j|^2.
\]
This implies
\[
 \{v\in B_i:|H(w)v|\le\tau\}
 \subseteq
 \left\{
 \sum_{j=1}^d z_j f_j:
 |z_j|\le
 \min\left\{2,\frac{\tau}{\Gamma_j(w)}\right\},
 \quad j=1,\ldots,d
 \right\},
\]
with $\tau/\Gamma_j(w)=+\infty$ when $\Gamma_j(w)=0$.
It follows that
\[
\begin{aligned}
 \left|\{v\in B_i:|H(w)v|\le\tau\}\right|
 &\lesssim
 \prod_{\ell=1}^{d}
 \min\left\{1,\frac{\tau}{\Gamma_\ell(w)}\right\}\\
 &=\tau^dD_\tau(w)^{-1}.
\end{aligned}
\]
Hence
\begin{equation}\label{eq:incidence-average-upper}
 F(w)\lesssim\tau^{dk}D_\tau(w)^{-k}.
\end{equation}

Whenever $F(w)>0$, there exist $v_i\in B_i$ satisfying
$|H(w)v_i|\le\tau$. Their uniform independence gives
\[
 \left|H(w)\sum_{i=1}^{k}a_i v_i\right|
 \le\tau\sum_{i=1}^{k}|a_i|
 \le2\sqrt{k}\,\tau
 \left|\sum_{i=1}^{k}a_i v_i\right|.
\]
Thus $H(w)$ has operator norm at most $2\sqrt{k}\,\tau$
on the $k$-dimensional space
$E=\operatorname{span}\{v_1,\ldots,v_k\}$.
By the min--max characterization of singular values,
\[
 \Gamma_i(w)\le2\sqrt{k}\,\tau,
 \qquad d-k+1\le i\le d.
\]
Since the remaining singular values are uniformly bounded,
\[
 D_\tau(w)\lesssim\tau^k
 \qquad\text{whenever }F(w)>0.
\]
For $s\ge k$, equation \eqref{eq:incidence-average-upper}
therefore yields
\[
 F(w)
 \lesssim
 \tau^{dk}D_\tau(w)^{s-k}D_\tau(w)^{-s}
 \lesssim
 \tau^{dk+k(s-k)}D_\tau(w)^{-s}.
\]
Combining this with \eqref{eq:incidence-average-lower},
we conclude that
\[
\begin{aligned}
 \int_{1<|w|<2}D_\tau(w)^{-s}\,dw
 &\gtrsim
 \tau^{N_k-dk-k(s-k)}\\
 &=\tau^{-ks+k(k+1)/2},
\end{aligned}
\]
as claimed.
\end{proof}
We now state a key lemma used throughout the proofs
to pass from smooth densities to general distributions
under the appropriate hypotheses.
\begin{lemma}[Annular testing]
\label{lem:annular-testing}
Let $1\le d<n$, put $k=n-d$, and let
\[
 M=\{(t',\varphi(t')):t'\in U\}\subset\mathbb R^d\times\mathbb R^k,
 \qquad \varphi\in C^\infty(U;\mathbb R^k),
\]
where $U\subset\mathbb R^d$ is open and $\sup_{U}\| D\phi \|\leq \frac{1}{10}$. Let $\kappa_\sigma$ denote the
annular Fourier-ratio decay exponent associated with $M$, and set
\[
 \Gamma_R=\{x\in\mathbb R^n:R<|x|<2R\}.
\]
Let $0\ne T\in\mathcal S'(\mathbb R^n)$ be supported in a compact
subset of $M$. Then there exists an
integer $j=j(T)\ge0$ such that, for every $2\le p<\infty$ and every
$\varepsilon>0$, there exist constants $c_{T,p,\varepsilon}>0$ and
$R_0(T,p,\varepsilon)\ge1$ for which
\begin{equation}
\label{eq:annular-lower}
 \|\hat{T}\|_{L^p(\Gamma_R)}
 \ge c_{T,p,\varepsilon}
 R^{\,j+\frac{n-2\kappa_\sigma}{p}
             -\frac{d-2\kappa_\sigma}{2}-\varepsilon},
 \qquad R\ge R_0(T,p,\varepsilon).
\end{equation}
\end{lemma}

\begin{proof}
Define the smooth flattening map
$S:U\times\R^k\rightarrow U\times\R^k$ by
\[
 S(t',t'')=(t',t''-\varphi(t')),
\]
and define its pushforward $T_s=S_*T$ by
\[
 \langle T_s,\Psi\rangle
 :=\langle T,\Psi\circ S\rangle
 ,
 \qquad \Psi\in C_c^\infty(\R^n).
\]
Since $\operatorname{supp}T_s\subset\{t''=0\}$, the structure
theorem for distributions supported on a subspace gives
\[
 T_s=\sum_{|\alpha|\le m}
 T_\alpha\otimes\partial^\alpha\delta_0,
 \qquad T_\alpha\in\mathcal E'(\R^d)
.\]
Choose $\alpha_0$ with $T_{\alpha_0}\ne0$ and $|\alpha_0|=m$.
We may choose $a\in C_c^\infty(U)$ with
$\langle T_{\alpha_0},a\rangle\ne0$ and by prescribing finitely
many moments of $\widehat b$, choose
$b\in\mathcal S(\R^k)$ such that
\[
 \operatorname{supp}\widehat b
 \Subset\{\lambda:4/3<|\lambda|<3/2\},
 \qquad
 \partial^\alpha b(0)=
 \begin{cases}
  1, & \alpha=\alpha_0,\\
  0, & \alpha\ne\alpha_0,\quad |\alpha|\le m.
 \end{cases}
\]
Define $g_R\in\mathcal S(\R^n)$ by
\[
 g_R(t',t'')=
 \begin{cases}
  a(t')\,b\bigl(R(t''-\varphi(t'))\bigr), & t'\in U,\\
  0 & t'\notin U.
 \end{cases}
\]
By definition of $T_s$,
\begin{equation}
\label{eq:curved-pairing}
\begin{aligned}
 \langle T,g_R\rangle
 &=\langle T_s,a(t')b(Rt'')\rangle\\
 &=\sum_{|\alpha|\le m}(-1)^{|\alpha|}R^{|\alpha|}
     \langle T_\alpha,a\rangle\,\partial^\alpha b(0)
 =cR^m,
\end{aligned}
\end{equation}
where $m=|\alpha_0|$ and
$c=(-1)^m\langle T_{\alpha_0},a\rangle\ne0$. First, we estimate $\widehat g_R$ and use Fourier duality to
relate this estimate to $\widehat T$.
A direct computation yields
\begin{equation}
\label{eq:test-transform}
 \widehat g_R(x)=R^{-k}\widehat b(x''/R)F_a(x),
 \qquad
 F_a(x)=\int_U
 e^{-i(t'\cdot x'+x''\cdot\varphi(t'))}a(t')\,dt'.
\end{equation}
Since $\widehat b(x''/R)$ restricts the transverse variable to
$4R/3<|x''|<3R/2$, we estimate $\widehat g_R$ on a suitable
annulus at scale $R$ and show that its tail outside this annulus
decays rapidly.

Choose a radial $\theta\in C_c^\infty(\R^n)$ supported in
$\{1\le|x|\le2\}$, with $0\le\theta\le1$ and
$\theta=1$ on $\{5/4\le|x|\le7/4\}$, and set
$\theta_R(x)=\theta(x/R)$.
We may then write
\[
 \widehat g_R
 =\theta_R\widehat g_R+(1-\theta_R)\widehat g_R.
\]
To estimate $\theta_R\widehat g_R$, we use
\eqref{eq:test-transform} followed by the Fourier-ratio bound
to obtain
\[
 \|\theta_R\widehat g_R\|_1
 \le R^{-k}\|\widehat b\|_\infty
 \|F_a\|_{L^1(\Gamma_R)}
 \lesssim_{a,b,\varepsilon}
 R^{\frac n2-k-\kappa_\sigma+\varepsilon}
 \|F_a\|_{L^2(\Gamma_R)}.
\]
We use the $L^2$ density bound\cite{AH76},
to control the $L^{2}$ norm of $F_{a}$ to get
\[
 \|\theta_R\widehat g_R\|_1
 \lesssim_{a,b,\varepsilon}
 R^{\frac d2-\kappa_\sigma+\varepsilon}.
\]
An application of Plancherel's theorem yields
\[
 \|\theta_R\widehat g_R\|_2
 \le \|\widehat g_R\|_2
 \lesssim_{a,b} R^{-k/2}.
\]
Combining the above two estimates by interpolation, with
$p'=p/(p-1)$, gives
\begin{equation}
\label{eq:dual-test-bound}
\begin{aligned}
 \|\theta_R\widehat g_R\|_{p'}
 &\le \|\theta_R\widehat g_R\|_1^{\,1-2/p}
       \|\theta_R\widehat g_R\|_2^{\,2/p}\\
 &\lesssim_{a,b,p,\varepsilon}
 R^{\frac{d-2\kappa_\sigma}{2}
       -\frac{n-2\kappa_\sigma}{p}+\varepsilon}.
\end{aligned}
\end{equation}
On the support of $(1-\theta_R)$, $\widehat b(x''/R)$
vanishes in the region $|x|\le5R/4$, so we estimate the decay
of $(1-\theta_R)\widehat g_R$ in the region $|x|\ge7R/4$.
In particular, we may use a nonstationary phase estimate
for $F_a$, since
\[
 \left|\nabla_{t'}
 \bigl(x'\cdot t'+x''\cdot\varphi(t')\bigr)\right|
 \ge |x'|-\tfrac1{10}|x''|
 \ge  c(|x|+R),
\]
where $c>0$ follows from the support restriction on
$\widehat b(x''/R)$ and the bound
$\sup_U\|D\varphi\|\le1/10$. All positive-order derivatives of the phase divided by
$R+|x|$ are uniformly bounded on the fixed support of $a$.
Nonstationary phase therefore gives, for every multi-index
$\alpha$ and every $L\ge1$,
\[
 |\partial_x^\alpha F_a(x)|
 \le C_{a,\alpha,L}(R+|x|)^{-L},
 \qquad
 \bigl(|x|\ge7R/4,\;
 x''/R\in\operatorname{supp}\widehat b\bigr),
\]
uniformly for $R\ge1$.

Since $\widehat T$ is tempered, there is an integer $q\ge0$
such that
\[
 |\langle \widehat T,F\rangle|
 \le C_T\sum_{|\alpha|\le q}
 \sup_x(1+|x|)^q|\partial^\alpha F(x)|,
 \qquad F\in\mathcal S(\R^n).
\]
Applying this to $F=(1-\theta_R)\widehat g_R(-\cdot)$,
Leibniz's rule and the preceding nonstationary phase estimate,
with $L>q+N$, give
\begin{equation}
\label{eq:test-tail}
 \left|\left\langle \widehat T,
 (1-\theta_R)\widehat g_R(-\cdot)\right\rangle\right|
 \lesssim_{T,a,b,N}R^{-N},
 \qquad R\ge1, N\geq0.
\end{equation}
To pass the estimates for $\widehat g_R$ to $\widehat T$,
we first use Fourier duality to obtain
\[
\begin{aligned}
 \langle T,g_R\rangle
 &=(2\pi)^{-n}
   \langle\widehat T,\widehat g_R(-\cdot)\rangle\\
 &=(2\pi)^{-n}
   \int_{\Gamma_R}\widehat T(x)\theta_R(x)
   \widehat g_R(-x)\,dx +(2\pi)^{-n}
   \left\langle\widehat T,
   (1-\theta_R)\widehat g_R(-\cdot)\right\rangle.
\end{aligned}
\]
Hölder's inequality followed by \eqref{eq:dual-test-bound}
yields
\[
 \left|\int_{\Gamma_R}\widehat T(x)\theta_R(x)
 \widehat g_R(-x)\,dx\right| \le C_{T,p,\varepsilon}
 R^{\frac{d-2\kappa_\sigma}{2}
       -\frac{n-2\kappa_\sigma}{p}+\varepsilon}
 \|\widehat T\|_{L^p(\Gamma_R)}.
\]
Combining this with \eqref{eq:test-tail}, we obtain
\[
 |\langle T,g_R\rangle|
 \lesssim_{T,p,\varepsilon} 
 R^{\frac{d-2\kappa_\sigma}{2}
       -\frac{n-2\kappa_\sigma}{p}+\varepsilon}
 \|\widehat T\|_{L^p(\Gamma_R)}
 +R^{-N}.
\]
It follows from \eqref{eq:curved-pairing},
$$
 R^j\lesssim_{T,p,\varepsilon}
 R^{\frac{d-2\kappa_\sigma}{2}
       -\frac{n-2\kappa_\sigma}{p}+\varepsilon}
 \|\hat{T}\|_{L^p(\Gamma_R)}+R^{-N}.
$$
For large $R$, absorb the remainder into $|c|R^m/2$ and rearrange to obtain
$$ \|\hat{T}\|_{L^p(\Gamma_R)}
 \gtrsim_{T,p,\varepsilon}
 R^{\,j+\frac{n-2\kappa_\sigma}{p}
             -\frac{d-2\kappa_\sigma}{2}-\varepsilon},
 \qquad R\ge R_0(T,p,\varepsilon).$$
\end{proof}
\section{Appendix}\label{Appendix}
We present a proof of Proposition\ref{Prop:Uniformrankdecay}
\begin{proof}
Let $\xi\in\mathbb R^d$, $\eta\in\mathbb R^{n-d}$, and set
$t=|(\xi,\eta)|$. We divide the proof into the
region $|\xi|\gg|\eta|$ and $|\eta| \gg |\xi|$, where constants are chosen sufficiently large depending upon $\phi$, and $A$. More precisely, set
\[
B=1+\sup_{x\in\operatorname{supp}(A)}\|D\Phi(x)\|.
\]
If $|\xi|>2B|\eta|$, then we have no critical points, thus we may repeat integration by parts $\lceil\frac{N}{2}\rceil$ times to get the desired bounds. Suppose now that $|\xi|\leq2B|\eta|$. Define
\[
\omega=\frac{1}{t}(\xi,\eta)
=(\omega',\omega'')\in\mathbb S^{n-1},
\qquad
g(\omega,x)=\langle\omega,(x,\Phi(x))\rangle.
\]
Since each $\phi_j\in C^{d+3}$, $g$ has the same regularity. The directions in this region belong to the compact  cone
\[
\Gamma=\{\omega\in\mathbb S^{n-1}:
|\eta| \gg |\xi|\}.
\]
Since $\eta\neq0$ on $\Gamma$, our hypothesis gives
\[
\operatorname{rank}D_x^2g(\omega,x)
=
\operatorname{rank}\left(
\sum_{j=1}^{n-d}\eta_{j}D^2\phi_j(x)\right)
\geq N
\]
for every $(\omega,x)\in\Gamma\times\operatorname{supp}(A)$.
The $N$th largest absolute eigenvalue of $D_x^2g(\omega,x)$ is
continuous in $(\omega,x)$. By compactness, there exists $\kappa>0$
such that
\[
|\lambda_N(\omega,x)|\geq\kappa
\qquad
\text{on }\Gamma\times\operatorname{supp}(A).
\]

Fix $(\omega_0,x_0)$ in this compact set. Let $Q_0$ be an orthogonal
matrix whose columns are eigenvectors of $D_x^2g(\omega_0,x_0)$,
ordered so that
\[
|\lambda_1|\geq\cdots\geq|\lambda_d|.
\]
Making the change of variables $x=x_0+Q_0y$, define
\[
g_0(\omega,y)=g(\omega,x_0+Q_0y).
\]
Then
\[
D_y^2g_0(\omega_0,0)
=
Q_0^TD_x^2g(\omega_0,x_0)Q_0
=
\operatorname{diag}(\lambda_1,\ldots,\lambda_d).
\]
Write $y=(u,v)\in\mathbb R^N\times\mathbb R^{d-N}$
By continuity, we may choose a neighborhood of $(\omega_0,0)$
on which
\[
|\det D_u^2g_0(\omega,u,v)|
\geq\left(\frac{\kappa}{2}\right)^N.
\]

By compactness, finitely many such neighborhoods cover
$\Gamma\times\operatorname{supp}(A)$. Choose a smooth and finite partition of
unity $\{\psi_i(\omega,x)\}_{i\in J}$ subordinate to this cover.
For the neighborhood centered at $(\omega_i,x_i)$, use the fixed
rotation $Q_i$ and write
\[
g_i(\omega,y)=g(\omega,x_i+Q_i y),
\qquad
A_i(\omega,y)
=\psi_i(\omega,x_i+Q_i y)A(x_i+Q_i y).
\]
Thus
\[
\int_{\mathbb R^d}e^{it g(\omega,x)}A(x)\,dx
=
\sum_{i\in J}
\int_{\mathbb R^{d-N}}\int_{\mathbb R^N}
e^{it g_i(\omega,u,v)}A_i(\omega,u,v)\,du\,dv.
\]
For a fixed $v \in \R^{d-N}$, applying stationary phase estimate for inner integral in $u$\cite{Hormander}, together with nonstationary phase away from
the critical points, gives
\[
\left|\int_{\mathbb R^N}
e^{it g_i(\omega,u,v)}A_i(\omega,u,v)\,du\right|
\lesssim_{A,\Phi}\kappa^{-N/2}t^{-N/2};
\]
Let $V_i$ be a fixed bounded set containing the $v$-projection of
the support of $A_i(\omega,\cdot)$ for all relevant $\omega$.
It follows that
\[
\begin{aligned}
\left|\int_{\mathbb R^{d-N}}\int_{\mathbb R^N}
e^{it g_i(\omega,u,v)}A_i(\omega,u,v)\,du\,dv\right|
&\leq
\int_{V_i}
\left|\int_{\mathbb R^N}
e^{it g_i(\omega,u,v)}A_i(\omega,u,v)\,du\right|\,dv\\
&\lesssim_{A,\Phi}t^{-N/2}.
\end{aligned}
\]
Summing over the finite cover and using
$t\sim1+|\xi|+|\eta|$ for $t>1$, we conclude that
\[
\left|\int_{\mathbb R^d}
e^{i(\xi\cdot x+\eta\cdot\Phi(x))}A(x)\,dx\right|
\lesssim_{A,\Phi}(1+|\xi|+|\eta|)^{-N/2}.
\]
\end{proof} 

\textbf{Acknowledgment and Declaration of AI use}
I would like to thank my advisor, Alex Iosevich(AI), for his helpful
suggestions and guidance throughout the preparation of this manuscript.
I would also like to thank Jon Pakianathan for helpful
observations regarding Hessians and related Linear Algebra.

AI tools were used primarily for grammar checks, polishing, TikZ
diagrams, and literature searches. The main results of this paper are
my own, though I also used AI tools to discuss methods, the scope of
the results, and alternative approaches to certain steps. In particular,
my original argument in Step~IV of the proof of
Theorem~\ref{thm:smooth-high-codimension} was substantially longer and
used a probabilistic method. Discussions with AI led to the identification
of exact optimizers, which simplified the argument and its presentation. I take full responsibility for the mathematical
claims, proofs, and references in this manuscript.

\enlargethispage{2\baselineskip}
\end{document}